\documentclass[11pt,draft]{amsart}

\newtheorem{theorem}{Theorem}[section]
\newtheorem{lemma}[theorem]{Lemma}

\theoremstyle{definition}
\newtheorem{defn}[theorem]{Definition}
\newtheorem{example}[theorem]{Example}

\newtheorem{remark}[theorem]{Remark}

\numberwithin{equation}{section}

\RequirePackage{amsmath}
\RequirePackage{amssymb}
\RequirePackage{amsthm}
 \RequirePackage{epsfig}

\begin{document}

\date{}

\title[A new approach to the real numbers]
{A new approach to the real numbers}

\author{Liangpan Li}

\address{Department of Mathematics, Shanghai Jiao Tong University,
Shanghai 200240,  China}

\address{Department of Mathematical Sciences, Loughborough
University, LE11 3TU, UK}
 \email{liliangpan@yahoo.com.cn}



\date{}

\begin{abstract}
In this paper we provide a complete approach to the real numbers via
decimal representations. Construction of the real numbers by
Dedekind cuts, Cauchy sequences of rational numbers, and the
algebraic characterization of the real number system by the  concept
of complete ordered field are also well explained in the new
setting.
\end{abstract}

\maketitle

\section{Dedekind cuts, Cauchy sequences, axiomatic approach and decimal representations}\label{section1}

``\textsf{The discovery of incommensurable quantities was a severe
blow to the Pythagorean program of understanding nature by means of
numbers. (The slogan of the Pythagoreans was `All is number.') The
Greeks developed a sophisticated theory of ratios, presumably the
work of Eudoxus, to work around the problem that certain quantities,
even certain lengths, could not be reduced to numbers, as then
understood. This theory anticipates the development of the real
number system by Dedekind and Cantor in the nineteenth century.}"
(\cite{Browder})

 In popular
literatures  there are mainly three approaches to the real numbers
such as  construction by Dedekind cuts (see e.g.
\cite{Browder,Burrill,Chen,Dieudonne2,Eves,Fikh,Hardy,Khinchin,Kodaria,Rana,Rudin}),
Cauchy sequences of rational numbers (see e.g.
\cite{Bridges,Burrill,Chen,Cohen,Dieudonne2,Hewitt,Rana,Roberts,Tao,Thurston}),
and an axiomatic definition (see e.g.
\cite{Apostol,Bridges,Browder,Burrill,Dieudonne,Douglass,Eves,Kirkwood,Korner,Lewin,Pedrick,Rana,Roberts,Royden,Rudin,Stoll,Zorich}).
We first make a brief summary.

Let $\mathbb{Q}$ be the set of rational numbers endowed with the
standard additive operation $+$, multiplicative operation $\times$,
also denoted by $\cdot$, and total order $\leq$. A \textsf{Dedekind
cut} is a pair of nonempty subsets $A,B$ of $\mathbb{Q}$, denoted by
$(A|B)$, such that
\begin{itemize}
\item $A\cap B=\emptyset$, $A\cup B=\mathbb{Q}$,
\item $a\in A, b\in B\Rightarrow a<b$,
\item $A$ contains no greatest element.
\end{itemize}
The set of all Dedekind cuts is denoted by $\mathbf{DR}$. In 1872,
Richard Dedekind (\cite{Dedekind,Dedekind2}) introduced an operation
$\boxplus$ called ``\textsf{addition}", an operation $\boxtimes$
called ``\textsf{multiplication}" and a total order $\lesssim$ on
$\mathbf{DR}$, getting the so-called real number system
$(\mathbf{DR},\boxplus,\boxtimes,\lesssim)$. This work profoundly
influenced subsequent studies of the foundations of mathematics.

A  sequence of rational numbers $\{x^{(n)}\}_{n\in\mathbb{N}}$ is
called \textsf{Cauchy} if for every  rational $\epsilon>0$, there
exists an $N\in\mathbb{N}$ such that  $\forall m,n \geq N$,
$|x^{(m)}-x^{(n)}|<\epsilon$. The set of all Cauchy sequences of
rational numbers is denoted by $\mathbf{CR}$.  Two Cauchy sequences
of rational numbers $\{x^{(n)}\}_{n\in\mathbb{N}}$,
$\{y^{(n)}\}_{n\in\mathbb{N}}$ are said to be equivalent, denoted by
$\{x^{(n)}\}_{n\in\mathbb{N}}\approx\{y^{(n)}\}_{n\in\mathbb{N}}$,
if for every rational $\epsilon>0$, there exists an $N\in\mathbb{N}$
such that $\forall n \geq N$, $|x^{(n)}-y^{(n)}|<\epsilon$.  It is
easy to prove $\approx$ is an equivalence relation, thus yields a
quotient space $\mathbf{CR}/\approx$ denoted by
$\widehat{\mathbf{CR}}$ for consistency. Around 1869$\sim$1872,
Charles M\'{e}ray (\cite{Meray}), Georg Cantor (\cite{Cantor}) and
Eduard Heine (\cite{Heine}) independently introduced basically the
same additive operation $\widehat{\boxplus}$, multiplicative
operation $\widehat{\boxtimes}$ and total order $\widehat{\lesssim}$
on $\widehat{\mathbf{CR}}$, getting yet another number system
$(\widehat{\mathbf{CR}},\widehat{\boxplus},\widehat{\boxtimes},\widehat{\lesssim})$.
This approach has the advantage of providing a standard way for
completing an abstract metric space.

Both systems are isomorphic in the sense that one can find a
bijection $\omega:\mathbf{DR}\rightarrow\widehat{\mathbf{CR}}$ such
that for any two elements $x,y\in\mathbf{DR}$, $\omega(x\boxplus
y)=\omega(x)\widehat{\boxplus}\ \omega(y)$, $\omega(x\boxtimes
y)=\omega(x)\widehat{\boxtimes}\ \omega(y),$ $x\lesssim
y\Leftrightarrow \omega(x)\widehat{\lesssim}\ \omega(y)$.
 Detailed studies of a dozen or so properties of $\boxplus,\boxtimes,\lesssim$ lead to the
basic concept of ``\textsf{complete ordered field}" and an
algebraic-axiomatic approach to the real number system. Later on,
when to verify a new system is isomorphic to those of Dedekind and
M\'{e}ray-Cantor-Heine, it suffices to prove that it is a complete
ordered field. We should also note several severe criticisms of the
algebraic-axiomatic approach:

``\textsf{The necessary axioms should come as a byproduct of the
construction process and not be predetermined.}" (\cite{Leviatan})

``\textsf{$\ldots$ the algebraic-axiomatic definition of a real
number is simply appalling and abhorrent $\ldots$ to define a real
number via a cold and boring list of a dozen or so axioms for a
complete ordered field is like replacing life by death or reading an
obituary column.}" (\cite{Abian})

$\sim\sim\sim\sim\sim\sim\sim\sim\sim\sim\sim\sim\sim\sim\sim\sim\sim\sim\sim\sim\sim\sim\sim\sim\sim\sim$

``\textsf{Few mathematical structures have undergone as many
revisions or have been presented in as many guises as the real
numbers. Every generation reexamines the reals in the light of its
values and mathematical objectives.}" (\cite{Faltin-Rota})

Next we introduce the real numbers through a rather old geometric
approach. Given a point $x$ in an ``axis", the decimal
representation of it is obtained step by step as follows.
\begin{itemize}
\item Step 0: Partition this ``axis" into countably many disjoint
unions $\bigcup_{z\in\mathbb{Z}}[z,z+1)$, then find a unique integer
$x_0\in\mathbb{Z}$ such that $x\in[x_0,x_0+1)$.
\item Step 1: Partition $[x_0,x_0+1)$ into ten disjoint unions
$\bigcup_{i=0}^9[x_0+\frac{i}{10},x_0+\frac{i+1}{10})$, then find a
unique element $x_{1}\in\mathcal{Z}_{10}$ such that
$x\in[x_0+\frac{x_1}{10},x_0+\frac{x_1+1}{10})$, here and afterwards
we denote $\{0,1,2,3,4,5,6,7,8,9\}$ by $\mathcal{Z}_{10}$ for
simplicity.
\item Step 2: Partition $[x_0+\frac{x_1}{10},x_0+\frac{x_1+1}{10})$ into ten
disjoint unions
$\bigcup_{j=0}^9[x_0+\frac{x_1}{10}+\frac{j}{10^2},x_0+\frac{x_1}{10}+\frac{j+1}{10^2})$,
then find a unique element $x_{2}\in\mathcal{Z}_{10}$ such that
$x\in[x_0+\frac{x_1}{10}+\frac{x_2}{10^2},x_0+\frac{x_1}{10}+\frac{x_2+1}{10^2})$.
\item Step $\vdots$: $\cdots\cdots$
\end{itemize}
Then we say $x$ has decimal representation $x_0.x_1x_2x_3\cdots$,
and call $x_k$ the $k$-th digit of $x$. A thorough consideration
leads to a natural question: can any point of this ``axis" have
decimal representation $0.999999\cdots$? and if not, why can we
expel its existence? This is not a silly question at all. To
correctly answer it, we should first know what an ``axis" it is, or
what on earth a ``real number" it is!

Although constructing  real numbers via decimals has been known
since Simon Stevin (\cite{Construction,Simon}) in the 16-th century,
developed also by Karl Weierstrass (\cite{Weierstrass}), Otto Stolz
(\cite{Stolz,Stolz2}) in the 19-th century, and many other modern
mathematicians (see e.g.
\cite{Abian,Behrend,Bruijn,Burrill,Courant,Faltin-Rota,Gamelin,Gowers,Hua,Klazar,Kudryavtsev,Leviatan,Lightstone,Richman,Wang}),
almost all popular Mathematical Analysis books didn't choose this
approach. The decimal construction hasn't got the attention it
should deserve:

 ``\textsf{Perhaps one of the greatest achievements
of the human intellect throughout the entire history of the human
civilization is the introduction of the decimal notation for the
purpose of recording the measurements of various magnitudes. For
that purpose the decimal notation is most practical, most simple,
and in addition, it reflects most outstandingly the profound
subtleties of the human analytic mind. In fact, decimal notation
reflects so much of the Arithmetic and so much of the Mathematical
Analysis $\ldots$}" (\cite{Abian})

To the author's opinion, at least one of the reasons behind this
phenomenon is most of the authors only gave outlines or sketches. We
note an interesting phenomenon happened in  popular Mathematical
Analysis books: many authors (see e.g.
\cite{Apostol,Browder,Chen,Fikh,Kalapodi,Kodaria,Korner,Pedrick,Rana,Roberts,Rudin,Stevenson,Stoll,Tao,Zorich})
first chose one of the other three approaches discussed before, then
proved that every real number has a suitable decimal representation.
But if without using the algebraic intuition that the real number
system is nothing but a complete ordered field which is not so easy
to grasp for beginners, we could find basically no literature
reversing this kind of discussion. No matter how fundamental
\textsf{sets} and \textsf{sequences} are in mathematics, conflicting
with our primary, secondary and high school education that a number
is a string of decimals, skepticism over explicit construction of
the real numbers by either Dedekind cuts or Cauchy sequences never
ends:

``\textsf{The degeometrization of the real numbers was not carried
out without skepticism. In his opus \emph{Mathematical Thought from
Ancient to Modern Times}, mathematics historian Morris Klein quotes
Hermann Hankel who wrote in 1867: Every attempt to treat the
irrational numbers formally and without the concept of [geometric]
magnitude must lead to the most abstruse ad troublesome
artificialities, which, even if they can be carried through complete
rigor, as we have every right to doubt, do not have a right
scientific value}." (\cite{Gamelin})

``\textsf{The definition of a real number as a Dedekind cut of
rational numbers, as well as a Cauchy sequence of rational numbers,
is cumbersome, impractical, and $\ldots$, inconsequential for the
development of the Calculus or the Real Analysis.}" (\cite{Abian})

Based on the fundamental concept of ``\textsf{order}" and its
derived operations, in this paper we will provide a complete
approach to the real numbers via decimals, and some of our ideas are
new to the existing literatures. Also in this new setting,
construction of the real numbers by Dedekind cuts, Cauchy sequences
of rational numbers, and the algebraic characterization of the real
number system by the  concept of complete ordered field can be well
explained. The general strategy of our approach is as follows.

The starting point is to choose \[{\mathcal R}\triangleq
\big\{x_0.x_1x_2x_3\cdots\big|\
x_{0}\in\mathbb{Z},x_k\in\mathcal{Z}_{10}, k\in\mathbb{N}\big\}.\]
as our ambient space, which was already discussed before. Once
adopted this decimal notation, we suggest you imagine it in mind as
a series $\sum_{k=0}^{\infty}\frac{x_i}{10^i}$. There are mainly two
reasons why we choose this notation, one is instead of discussing
``\textsf{subtraction}", we shall focus on ``\textsf{additive
inverse}" which will be introduced in an elegant manner, the other
comes from the next paragraph.

Since as sets $\mathcal R$ is basically the same as
$\mathbb{Z}\times\mathcal{Z}_{10}^{\mathbb{N}}$, we can introduce a
lexicographical order $\preceq$ on $\mathcal R$, then prove
\textsf{the least upper bound property} and \textsf{the greatest
lower bound property} for $(\mathcal R, \preceq)$ from the same
properties for $(\mathbb{Z},\leq)$ as soon as possible. As
experienced readers should know, this would mean that $(\mathcal R,
\preceq)$ is ``complete". So in this complete setting, we can derive
five basic operations such as the supremum operation $\sup(\cdot)$,
the infimum operation $\inf(\cdot)$, the upper limit operation
$\overline{\mbox{LIMIT}}(\cdot)$, the lower limit operation
$\underline{\mbox{LIMIT}}(\cdot)$ and the limit operation
$\mbox{LIMIT}(\cdot)$.

It is very natural to define
\[x\oplus y\triangleq\mbox{LIMIT}\big(\{[x]_k+[y]_k\}_{k\in\mathbb{N}}\big)\]
for any two elements  $x,y\in\mathcal R$, where $[x]_k\triangleq
x_0.x_1x_2\cdots
x_k=\frac{\sum_{i=0}^kx_i\cdot10^{k-i}}{10^k}\in\mathbb{Q}$ is the
truncation of $x$ up to the $k$-th digit, so is $[y]_k$. As for the
definition of $[x]_k+[y]_k$, even a primary school student may know
how to do it, that is, for example,
\begin{align*}
&(-15).3456\\
+&\underline{(-18).6789}\\
&(-32).0245
\end{align*}

To define a multiplicative operation in a succinct way we need some
preparation. First we introduce  a signal map $\mbox{sign}:{\mathcal
R}\rightarrow\{0,1\}$ by
\begin{equation*} \mbox{sign}(x)\triangleq
\begin{cases}
\ 0 \ \ \mbox{if}\  x\succeq0.000000\cdots, \\
\ 1 \ \ \mbox{if}\  x\preceq(-1).999999\cdots.
\end{cases}
\end{equation*}
This map partitions $\mathcal R$ into two parts, one is
$\mbox{sign}^{-1}(0)$, understood as the positive part of $\mathcal
R$, the other is $\mbox{sign}^{-1}(1)$, understood as the negative
part of $\mathcal R$. Both parts are closed connected through an
``additive inverse" map
\[\Psi(x_0.x_1x_2x_3\cdots)\triangleq(-1-x_0).(9-x_1)(9-x_2)(9-x_3)\cdots,\]
which turns a positive element into a negative one, and vice versa.
The absolute value of $x$, denoted by $\|x\|$, is defined to be the
maximum over $x$ and $\Psi(x)$. Because of the wonderful formula
$x=\Psi^{({\mbox{sign}}(x))}(\|x\|), $ the author likes to call them
three golden flowers. Now for any two elements $x,y\in\mathcal R$,
we define their multiplication by
\[x\otimes y\triangleq\Psi^{\big(\mbox{sign}(x)+\mbox{sign}(y)\big)}\Big(\mbox{LIMIT}\big(\{[\|x\|]_k\cdot[\|y\|]_k\}_{k\in\mathbb{N}}\big)\Big).\]
In primary school we have already learnt how to define
$[\|x\|]_k\cdot[\|y\|]_k$.

We remark that motivated by the above definitions of addition and
multiplication on $\mathcal R$, similar operations will be
introduced on $\mathbf{DR}$ in a highly consistent way.

\textsf{At this stage, we have introduced a rough system $(\mathcal
R,\preceq,\oplus,\otimes)$ without any pain}. But unfortunately,
several well-known properties generally a standard addition and a
standard multiplication should have don't hold for $\oplus$ and
$\otimes$. To overcome these difficulties,  we will introduce an
equivalence relation $\sim$
 identifying $0.999999\cdots$ with $1.000000\cdots$, and the same
 like. \textsf{No matter adopting decimal,  binary or hexadecimal notation, no matter introducing such relations earlier or later, we
 cannot avoid doing it}.

With these preparations, we then verify the commutative, associate
and distribute laws, the existences of additive (multiplicative)
unit and inverse, and so on. \textsf{Almost all the verification
work depend only on a pleasant Lemma \ref{eqvi test improved}}.
Finally we define the set $\mathbb{R}$ of \textsf{real numbers} to
be the set of equivalent classes $\mathcal R/\sim$ with derived
operations $\widehat{\oplus}$, $\widehat{\otimes}$,
$\widehat{\preceq}$ from $\oplus$, $\otimes$, $\preceq$
respectively, thus yields our desired number system
$(\mathbb{R},\widehat{\oplus},\widehat{\otimes},\widehat{\preceq})$.

Now in the new setting $(\mathcal R, \preceq)$ with derived
operations such as $\sup(\cdot)$, $\inf(\cdot)$ for
\textsf{subsets}, and $\overline{\mbox{LIMIT}}(\cdot)$ and
$\underline{\mbox{LIMIT}}(\cdot)$ for \textsf{sequences} of
$\mathcal R$, we can further explain the construction of the real
numbers by Dedekind cuts and Cauchy sequences of rational numbers.
Given a Dedekind cut $(A|B)$, we obviously have $\sup A\preceq\inf
B$. It would be very nice if $\sup A\sim\inf B$, thus one can derive
a map $\tau$ from $\mathbf{DR}$ to $\mathbb{R}$ by sending $(A|B)$
to $[\sup A]=[\inf B]$. Later on we shall prove that this is indeed
the case. Given a Cauchy   sequence of rational numbers
$\{x^{(n)}\}_{n\in\mathbb{N}}$, we obviously have
$\underline{\mbox{LIMIT}}(\{x^{(n)}\}_{n\in\mathbb{N}})\preceq\overline{\mbox{LIMIT}}(\{x^{(n)}\}_{n\in\mathbb{N}})$.
It would be great if
$\underline{\mbox{LIMIT}}(\{x^{(n)}\}_{n\in\mathbb{N}})\sim\overline{\mbox{LIMIT}}(\{x^{(n)}\}_{n\in\mathbb{N}})$,
thus one can derive a map $\kappa$ from $\mathbf{CR}$ to
$\mathbb{R}$ by sending $\{x^{(n)}\}_{n\in\mathbb{N}}$ to
$[\underline{\mbox{LIMIT}}(\{x^{(n)}\}_{n\in\mathbb{N}})]=[\overline{\mbox{LIMIT}}(\{x^{(n)}\}_{n\in\mathbb{N}})]$.
Later on we shall prove that this  is also indeed the case.
Motivated by these observations, we will continue to prove that our
number system is isomorphic to those of Dedekind and
M\'{e}ray-Cantor-Heine.

We can explain Cauchy sequences of rational numbers in another way.
Given a Cauchy   sequence of rational numbers
$\{x^{(n)}\}_{n\in\mathbb{N}}$, if it could represent a ``real
number", then we have no doubt that any subsequence of it, say for
example a monotonically increasing or a monotonically decreasing
subsequence, should represent the same ``real number". Then we put
such a monotone subsequence of rational numbers in $(\mathcal R,
\preceq)$ whose existence is taken for granted at this time, from
either the least upper bound property or the greatest lower bound
property in the new setting, the interested readers definitely know
which decimal representation should be understood as the ``real
number" the original sequence represents. Simply speaking, it would
be great if we can discuss Dedekind cuts or Cauchy sequences of
rational numbers in a constructive, complete
 setting.

Still in the setting $(\mathcal R, \preceq)$, we will give the
traditional characterization of \textsf{irrational numbers}, which
makes our approach matching what we have learnt in high school. We
also explain how can we come to the basic concept of complete
ordered field.

$\sim\sim\sim\sim\sim\sim\sim\sim\sim\sim\sim\sim\sim\sim\sim\sim\sim\sim\sim\sim\sim\sim\sim\sim\sim\sim$

To develop this approach, the author owed a lot to Loo-Keng Hua's
masterpiece \cite{Hua}. He also thanks Rong Ma for helpful
discussions. This work was partially supported by the Natural
Science Foundation of China (Grant Number 11001174).

Some notations used throughout this paper:
\begin{itemize}
\item $\mathbb{N}$ is the set of natural numbers
\item $\mathbb{Z}$ is the set of integers
\item Denote $\{0,1,2,3,4,5,6,7,8,9\}$ by $\mathcal{Z}_{10}$

\item $\lfloor\cdot\rfloor$ is the floor function, $\lceil\cdot\rceil$
is the ceil function

\item For any nonempty bounded above (below) subset
$E$ of $\mathbb{Z}$, denote by $\max E$ ($\inf E$) the unique least
upper (greatest lower) bound for $E$

\item
For any binary operation $\circledast$ on a set $G$, denote by
$A\circledast B$ the set $\{a\circledast b|\ a\in A, b\in B\}$,
where $A,B$ are nonempty subsets of $G$

\item Rational number system: $(\mathbb{Q},+,\times,\leq)$

\item Decimal system: $(\mathcal
R,\oplus,\otimes,\preceq)$

\item Dedekind cut system: $(\mathbf{DR},\boxplus,\boxtimes,\lesssim)$

\item Cauchy sequence system: $(\widehat{\mathbf{CR}},\widehat{\boxplus},\widehat{\boxtimes},\widehat{\lesssim})$

\item Real number system: $(\mathbb{R},\widehat{\oplus},\widehat{\otimes},\widehat{\preceq})$
\end{itemize}

\section{Least upper bound and greatest lower bound properties for $(\mathcal R, \preceq)$}

\subsection{Ambient space, lexicographical order and its derived operations}
\begin{defn}
As in Section \ref{section1} we define the ambient space of this
paper to be
\[{\mathcal R}\triangleq \big\{x_0.x_1x_2x_3\cdots\big|\
x_{0}\in\mathbb{Z},x_k\in\mathcal{Z}_{10}, k\in\mathbb{N}\big\}.\]
For any $x=x_0.x_1x_2x_3\cdots\in{\mathcal R}$ and $k\in\mathbb{N}$,
let $[x]_k\triangleq x_0.x_1\cdots
x_k=\frac{\sum_{i=0}^kx_i\cdot10^{k-i}}{10^k}\in\mathbb{Q}$ be the
truncation of $x$ up to the $k$-th digit.
\end{defn}

\begin{remark}
For any  $\alpha\in\mathbb{Q}$, there exist two integers
$p\in\mathbb{N}$ and $q\in\mathbb{Z}$
 such that $\alpha=\frac{q}{p}$.  Via the \textsf{long
division algorithm} $\alpha$ has a decimal representation
$\alpha_0.\alpha_1\alpha_2\alpha_3\cdots$, where
\[\alpha_0\triangleq\lfloor\frac{q}{p}\rfloor,\ \beta_0\triangleq
q-\alpha_0p,\
\alpha_{k+1}\triangleq\lfloor\frac{10\beta_k}{p}\rfloor,\
\beta_{k+1}\triangleq10\beta_k-\alpha_{k+1}p.\] Obviously, this
representation is independent of the choices of $p$ and $q$, so from
now on we can always view $\mathbb{Q}$ as a subset of ${\mathcal
R}$. For example,
\[  \frac{5}{2}=2.5000000\cdots,\ \
\frac{-40}{3}=(-14).666666\cdots.\]
\end{remark}

\begin{defn}
Let $\preceq$ be the lexicographical order on $\mathcal R$, that is,
$x\preceq y$ if and only if $\forall k\in\mathbb{N}$,
$[x]_k\leq[y]_k$.
\end{defn}

As usual, we may write $y\succeq x$ when $x\preceq y$, and $x\prec
y$ when $x\preceq y$, $x\neq y$. Obviously, as elements of
$\mathbb{Q}$, $x\leq y$ if and only if $x\preceq y$, which is
self-evident since we can write $x,y$ by a common denumerator, then
use the long division algorithm.

\begin{defn}
A nonempty subset $W$ of ${\mathcal R}$ is called bounded above
(below) if there exists an $M\in{\mathcal R}$ such that $\forall
w\in W$, $w\preceq M$ $(w\succeq M)$. A nonempty subset of
${\mathcal R}$ is called bounded if it is bounded both above and
below. A sequence $\{x^{(n)}\}_{n\in\mathbb{N}}$ of $\mathcal R$ is
called  bounded above (below), bounded if the set $\{x^{(n)}|\
n\in\mathbb{N}\}$ is of the corresponding property.
\end{defn}

\begin{theorem}\label{main result}
Every nonempty bounded above subset of $({\mathcal R},\preceq)$ has
a least upper  bound, every nonempty bounded below subset of
$({\mathcal R},\preceq)$ has a greatest lower  bound.
\end{theorem}

\begin{proof}
This theorem follows from the least upper bound property and the
greatest lower bound property for $(\mathbb{Z},\leq)$. We shall only
prove the first part of this theorem, and leave the second one to
the interested readers. Let $W$ be a nonempty bounded above subset
of ${\mathcal R}$. Denote
\[M_{0}\triangleq\max\{x_0|\ x_0.x_1x_2x_3\cdots\in W\},\]
\[M_{k}\triangleq\max\{x_k|\ x_0.x_1x_2x_3\cdots\in W \  \mbox{with}\  x_i=M_i, i=0,1,\ldots,k-1\}\ \ \ (\forall k\in\mathbb{N}).\]
Obviously, $M\triangleq M_0.M_1M_2M_3\cdots$ is an upper bound for
$W$. Let $\widetilde{M}$ be an arbitrary upper bound for $W$.
$\forall k\in\mathbb{N}$, by the definition of $M_{k}$ one can find
an $x=x_0.x_1x_2x_3\cdots\in W$ such that $x_i=M_i$ for
$i=0,1,\ldots,k$. Thus $[\widetilde{M}]_k\geq[x]_k=x_0.x_1\cdots
x_k=M_0.M_1\cdots M_k=[M]_k$, which means $\widetilde{M}\succeq M$.
This proves $M$ is the (unique) least upper bound for $W$.

\end{proof}

\begin{defn}
The least upper bound, also known as the \textsf{supremum}, for a
nonempty bounded above subset $W$ of $\mathcal{R}$ is denoted by
$\sup W$, while the greatest lower bound, known as the
\textsf{infimum}, for a nonempty bounded below subset
$\widetilde{W}$ of $\mathcal{R}$ is denoted by $\inf \widetilde{W}$.
\end{defn}

\begin{defn}
Generally given a sequence $\{x^{(n)}\}_{n\in\mathbb{N}}$, we
should pay attention to its asymptotic behavior. Obviously for any
$n\in\mathbb{N}$, $\sup\{x^{(k)}|\ k\geq n\}$ can be understood as
an ``upper bound" for the sequence $\{x^{(n)}\}_{n\in\mathbb{N}}$ if
we don't care about several of its initial terms. Thus if we want to
obtain an asymptotic ``upper bound" that is as least as possible,
then we are naturally led to the concept of \textsf{upper limit} of
a \textsf{bounded} sequence, that is,
\[\overline{\mbox{LIMIT}}(\{x^{(n)}\}_{n\in\mathbb{N}})\triangleq
\inf\big\{\sup\{x^{(k)}|\ k\geq n\}\big|\ n\in\mathbb{N}\big\}.\]
Similarly, we define the \textsf{lower limit} of a \textsf{bounded}
sequence $\{x^{(n)}\}_{n\in\mathbb{N}}$ to be
\[\underline{\mbox{LIMIT}}(\{x^{(n)}\}_{n\in\mathbb{N}})\triangleq
\sup\big\{\inf\{x^{(k)}|\ k\geq n\}\big|\ n\in\mathbb{N}\big\}.\]
When it happens that
$\overline{\mbox{LIMIT}}(\{x^{(n)}\}_{n\in\mathbb{N}})=\underline{\mbox{LIMIT}}(\{x^{(n)}\}_{n\in\mathbb{N}})=L$,
we shall simply write $\mbox{LIMIT}(\{x^{(n)}\}_{n\in\mathbb{N}})$
for their common value $L$, and say the sequence
$\{x^{(n)}\}_{n\in\mathbb{N}}$ has limit $L$.
\end{defn}

\begin{remark}
It is no hard to prove that any bounded above monotonically
increasing sequence, that is, $x^{(1)}\preceq x^{(2)}\preceq
x^{(3)}\preceq x^{(4)}\preceq\cdots\preceq M$, or any bounded below
monotonically decreasing sequence, that is, $x^{(1)}\succeq
x^{(2)}\succeq x^{(3)}\succeq x^{(4)}\succeq\cdots\succeq M$, has a
limit. For example, suppose $x^{(1)}\preceq x^{(2)}\preceq
x^{(3)}\preceq x^{(4)}\preceq\cdots\preceq M$.  Then
\begin{align*}
\overline{\mbox{LIMIT}}(\{x^{(n)}\}_{n\in\mathbb{N}})&=
\inf\big\{\sup\{x^{(k)}|\ k\geq n\}\big|\ n\in\mathbb{N}\big\}\\
&=\inf\big\{\sup\{x^{(k)}|\ k\in\mathbb{N}\}\big|\
n\in\mathbb{N}\big\}\\
&=\sup\{x^{(k)}|\ k\in\mathbb{N}\}\\
&=\sup\big\{\inf\{x^{(n)}|\ n\geq k\}\big|\ k\in\mathbb{N}\big\}\\
&=\underline{\mbox{LIMIT}}(\{x^{(n)}\}_{n\in\mathbb{N}}).
\end{align*}
\end{remark}

We state below some elementary properties of $\sup(\cdot)$,
$\inf(\cdot)$, $\overline{\mbox{LIMIT}}(\cdot)$,
$\underline{\mbox{LIMIT}}(\cdot)$ and $\mbox{LIMIT}(\cdot)$ in a
lemma, and leave the proofs to the interested readers.

\begin{lemma}\label{lemma order}
\begin{enumerate}
\item $\inf W\preceq \sup W$,
\item $W_1\subset W_2\Rightarrow \sup W_1\preceq \sup W_2$,
\item $W_1\subset W_2\Rightarrow \inf W_1\succeq \inf W_2$,
\item
$\underline{\mbox{LIMIT}}(\{x^{(n)}\}_{n\in\mathbb{N}})\preceq\overline{\mbox{LIMIT}}(\{x^{(n)}\}_{n\in\mathbb{N}})$,
\item $x^{(n)}\preceq y^{(n)}\Rightarrow \overline{\mbox{LIMIT}}(\{x^{(n)}\}_{n\in\mathbb{N}})\preceq
\overline{\mbox{LIMIT}}(\{y^{(n)}\}_{n\in\mathbb{N}})$,
\item $x^{(n)}\preceq y^{(n)}\Rightarrow \underline{\mbox{LIMIT}}(\{x^{(n)}\}_{n\in\mathbb{N}})
\preceq\underline{\mbox{LIMIT}}(\{y^{(n)}\}_{n\in\mathbb{N}})$,
\item
$\overline{\mbox{LIMIT}}\big(\{x^{(n_i)}\}_{i\in\mathbb{N}}\big)\preceq\overline{\mbox{LIMIT}}\big(\{x^{(n)}\}_{n\in\mathbb{N}}\big)$,
\item
$\underline{\mbox{LIMIT}}\big(\{x^{(n_i)}\}_{i\in\mathbb{N}}\big)\succeq\underline{\mbox{LIMIT}}\big(\{x^{(n)}\}_{n\in\mathbb{N}}\big)$,
\item
$\mbox{LIMIT}\big(\{x^{(n)}\}_{n\in\mathbb{N}}\big)=L\Rightarrow\mbox{LIMIT}\big(\{x^{(n_i)}\}_{i\in\mathbb{N}}\big)=L$,
\item $\mbox{LIMIT}\big(\{[x]_k\}_{k\in\mathbb{N}}\big)=x$.
\end{enumerate}
\end{lemma}

\subsection{An equivalence relation}

\begin{defn}
Two elements $x,y$ of ${\mathcal R}$ are said to have no gap,
denoted by $x\sim y$, if it does not exist an element $z\in{\mathcal
R}\backslash\{x,y\}$ lying exactly between $x$ and $y$. \end{defn}

\begin{defn}
Let ${\mathcal R}_9$, $\mathbb{Q}_{F}$ be respectively the sets of
all decimal representations ending in an infinite string of nines,
and zeros.
\end{defn}

\begin{lemma}\label{R9Qf}
$x\prec y$, $x\sim y$ $\Rightarrow x\in {\mathcal R}_9, y\in
\mathbb{Q}_{F}$.
\end{lemma}
\begin{proof}
Suppose $x=x_0.x_1x_2x_3\cdots\prec y=y_0.y_1y_2y_3\cdots$. Let $k$
be the minimal non-negative integer such that $x_k<y_k$. Note
\[x\preceq x_0.x_1\cdots x_k999999\cdots\prec
y_0.y_1\cdots y_k000000\cdots\preceq y.\] Since $x\sim y$,  we must
have $x=x_0.x_1\cdots x_k999999\cdots$, else we get a contradiction
$x\prec x_0.x_1\cdots x_k999999\cdots\prec y$. Similarly,
$y=y_0.y_1\cdots y_k000000\cdots$. This finishes the proof.
\end{proof}

\begin{remark}
As a corollary of Lemma \ref{R9Qf},  it is easy to prove that $\sim$
is an equivalence relation on ${\mathcal R}$, and every equivalent
class has at most two elements of ${\mathcal R}$. As we know in real
life, one can find a medium of any two different points in a
straight line, but at this time it is hard for us to define the
medium between two points having no gap. So to express a straight
line from $\mathcal R$, it is absolutely necessary to module
something from $\mathcal R$, which we shall discuss in detail at a
very late stage of this paper.
\end{remark}

Next we prepare a useful characterization lemma. An enhanced Lemma
\ref{eqvi test improved} will be given in the next section.

\begin{lemma}\label{equivalence test}
 $x\sim y\Leftrightarrow\forall
k\in\mathbb{N}$, $|[x]_k-[y]_k|\leq\frac{1}{10^k}$.
\end{lemma}

\begin{proof}
The necessary part follows immediately from Lemma \ref{R9Qf}, so we
need only prove the sufficient one, and suppose $\forall
k\in\mathbb{N}$, $|[x]_k-[y]_k|\leq\frac{1}{10^k}$. Without loss of
generality we may assume that $x\preceq y$. Thus $\forall
k\in\mathbb{N}$, $[x]_k\leq[y]_k\leq[x]_k+\frac{1}{10^k}$.

\textsc{Case 1}: Suppose $x\in{\mathcal R}_9$. There exists an
$m\in\mathbb{N}$ such that $\forall k> m, x_k=9$. Obviously, $x\sim
x_0.x_1x_2\cdots x_m+\frac{1}{10^m}$. Note $\forall k> m$,
\[[x]_k\leq[y]_k\leq[x]_k+\frac{1}{10^k}=x_0.x_1x_2\cdots x_m+\frac{1}{10^m}.\]
Letting $k\rightarrow\infty$ gives $x\preceq y\preceq
x_0.x_1x_2\cdots x_m+\frac{1}{10^m}$, here we have used the fifth
and tenth parts of Lemma \ref{lemma order}. This naturally implies
$x\sim y$.

\textsc{Case 2}: Suppose $x\not\in{\mathcal R}_9$. There exists a
sequence of natural numbers $m_1<m_2<m_3<\cdots$ such that $\forall
i\in \mathbb{N}$, $x_{m_i}<9$. Noting
 $$[x]_{m_i}\leq[y]_{m_i}\leq[x]_{m_i}+\frac{1}{10^{m_i}}$$ and
 $$\big[[x]_{m_i}\big]_{m_i-1}=\big[[x]_{m_i}+\frac{1}{10^{m_i}}\big]_{m_i-1}=[x]_{m_i-1},$$
 we must have
$[y]_{m_i-1}=\big[[y]_{m_i}\big]_{m_i-1}=[x]_{m_i-1}.$ By the ninth
and tenth parts of Lemma \ref{lemma order}, we have $x=y$.

This concludes the whole proof of the lemma.
\end{proof}

\section{Additive operations}

\subsection{Additive operations}

\begin{defn}[Addition] For any two elements $x,y\in{\mathcal R}$,
let
\[x\oplus y\triangleq\mbox{LIMIT}\big(\{[x]_k+[y]_k\}_{k\in\mathbb{N}}\big).\]
This is well-defined since the sequence
$\{[x]_k+[y]_k\}_{k\in\mathbb{N}}$ is monotonically increasing with
an upper bound $x_0+y_0+2$, where $x=x_0.x_1x_2x_3\cdots$,
$y=y_0.y_1y_2y_3\cdots$.
\end{defn}

\begin{example}\label{example32}
For any $x\in{\mathcal R}$, $x\oplus0.000000\cdots=x$. This means
$({\mathcal R},\oplus)$ has a unit.
\end{example}

\begin{remark}
Given two elements $x,y\in\mathbb{Q}_F\subset\mathbb{Q}$, it is easy
to verify that $x\oplus y=x+y$. Therefore  from now on we can abuse
the uses of $\oplus$ and $+$ if the summands lie in $\mathbb{Q}_F$.
\end{remark}

\begin{theorem}\label{abel sum} For any two elements $x,y\in{\mathcal R}$, we have
$x\oplus y=y\oplus x.$
\end{theorem}

\begin{proof}
$x\oplus
y=\mbox{LIMIT}\big(\{[x]_k+[y]_k\}_{k\in\mathbb{N}}\big)=\mbox{LIMIT}\big(\{[y]_k+[x]_k\}_{k\in\mathbb{N}}\big)=y\oplus
x.$
\end{proof}

\begin{theorem}\label{keep order} Given $x,y,z,w\in{\mathcal R}$ with $x\preceq z$ and $y\preceq w$, we have
$x\oplus y\preceq z\oplus w.$
\end{theorem}

\begin{proof}
$\forall k\in\mathbb{N}$ we have $[x]_k\leq[z]_k$ and
$[y]_k\leq[w]_k$, which yields $[x]_k+[y]_k\leq[z]_k+[w]_k$. By
Lemma \ref{lemma order},
$\mbox{LIMIT}\big(\{[x]_k+[y]_k\}_{k\in\mathbb{N}}\big)\preceq\mbox{LIMIT}\big(\{[z]_k+[w]_k\}_{k\in\mathbb{N}}\big)$.
\end{proof}

\begin{lemma}\label{additive estimate one element KEY}
For any element $x\in{\mathcal R}$ and $k\in\mathbb{N}$, we have
$[x]_k\preceq x\preceq[x]_k+\frac{1}{10^k}$.
\end{lemma}

\begin{proof}The first inequality is evident and we need only to prove the
second one. For any natural numbers $n\geq k$,
$[x]_n\leq[x]_k+\frac{1}{10^k}$. Letting $n\rightarrow\infty$ yields
$x\preceq[x]_k+\frac{1}{10^k}$.
\end{proof}

\begin{lemma}\label{additive estimate two elements KEY}
For any elements $x,y\in{\mathcal R}$ and $k\in\mathbb{N}$, we have
$$[x]_k+[y]_k\leq [x\oplus y]_k\preceq x\oplus
y\preceq[x]_k+[y]_k+\frac{2}{10^k}.$$
\end{lemma}

\begin{proof}
To prove the first inequality, we need only note
$$[x]_k+[y]_k\preceq x\oplus
y\Rightarrow[x]_k+[y]_k=[[x]_k+[y]_k]_k\leq[x\oplus y]_k.$$ For any
natural numbers $n\geq k$,
\[[x]_n+[y]_n\leq([x]_k+\frac{1}{10^k})+([y]_k+\frac{1}{10^k})=[x]_k+[y]_k+\frac{2}{10^k}.\]
Letting $n\rightarrow\infty$ yields $x\oplus
y\preceq[x]_k+[y]_k+\frac{2}{10^k}$. This proves the third
inequality. The second one follows from the previous lemma, so we
finishes the whole proof.
\end{proof}

 \begin{lemma}
\label{eqvi test improved} If there exists an $M\in \mathbb{N}$ such
that  $\forall k\in\mathbb{N}$, $|[x]_k-[y]_k|\leq\frac{M}{10^k}$,
then $x\sim y$.
\end{lemma}

\begin{proof}
For any $k\in\mathbb{N}$,
$y\preceq[y]_k+\frac{1}{10^k}\leq[x]_k+\frac{M+1}{10^k}\preceq
x\oplus\frac{M+1}{10^k}$. For any $m\in \mathbb{N}$, we can find a
sufficiently large $k$ such that
$\frac{M+1}{10^k}\leq\frac{1}{10^m}$. Consequently, $y\preceq
x\oplus\frac{1}{10^m}$, which yields $[y]_m\leq
[x\oplus\frac{1}{10^m}]_m=[x]_m+\frac{1}{10^m}$. By symmetry we can
also have $[x]_m\leq[y]_m+\frac{1}{10^m}$. Thus
$|[x]_m-[y]_m|\leq\frac{1}{10^m}$, and this finishes the proof
simply by applying Lemma \ref{equivalence test}.
\end{proof}

\begin{theorem}\label{keep order 33} Given $x,y,z,w\in{\mathcal R}$ with $x\sim z$ and $y\sim w$, we have
$x\oplus y\sim z\oplus w.$
\end{theorem}

\begin{proof}
For any $k\in\mathbb{N}$,
\[[x\oplus
y]_k\leq[x]_k+[y]_k+\frac{2}{10^k}\leq[z]_k+[w]_k+\frac{4}{10^k}\leq[z\oplus
w]_k+\frac{4}{10^k},\] here we have used Lemma \ref{equivalence
test} for $x\sim z$ and $y\sim w$. By symmetry we can also have
$[z\oplus w]_k\leq[x\oplus y]_k+\frac{4}{10^k}.$ Finally by Lemma
\ref{eqvi test improved}, we are done.
\end{proof}

\begin{theorem}\label{zzzzz} For any three elements $x,y,z\in{\mathcal R}$, we have
$(x\oplus y)\oplus z\sim x\oplus(y\oplus z).$
\end{theorem}

\begin{proof}
For any $k\in\mathbb{N}$,
\begin{align*}[x]_k+[y]_k+[z]_k&\leq[x\oplus y]_k+[z]_k\leq[(x\oplus y)\oplus z]_k\preceq(x\oplus y)\oplus
z\\&\preceq([x]_k+[y]_k+\frac{2}{10^k})+([z]_k+\frac{1}{10^k})=[x]_k+[y]_k+[z]_k+\frac{3}{10^k}.
\end{align*}
Similarly, we can also have $[x]_k+[y]_k+[z]_k\leq[(y\oplus z)\oplus
x]_k\leq[x]_k+[y]_k+[z]_k+\frac{3}{10^k}.$ Finally by Lemma
\ref{eqvi test improved}, $(x\oplus y)\oplus z\sim(y\oplus z)\oplus
x=x\oplus(y\oplus z)$. This concludes the proof.
\end{proof}

\begin{remark}
Given $n$ elements $\{x^{(i)}\}_{i=1}^n$ of ${\mathcal R}$ and a
permutation $\tau$ on the index set $\{1,2,\ldots,n\}$, according to
Theorems \ref{abel sum}, \ref{keep order 33} and \ref{zzzzz}, it is
easy to verify that
$$(\cdots(((x^{(1)}\oplus x^{(2)})\oplus x^{(3)})\oplus x^{(4)})\cdots)\oplus x^{(n)}\sim
(\cdots(((x^{(\tau_1)}\oplus x^{(\tau_2)})\oplus x^{(\tau_3)})\oplus
x^{(\tau_4)})\cdots)\oplus x^{(\tau_n)}.$$ As usual, we may simply
write
$$x^{(1)}\oplus x^{(2)}\oplus x^{(3)}\oplus x^{(4)}\oplus \cdots\oplus x^{(n)}\sim
x^{(\tau_1)}\oplus x^{(\tau_2)}\oplus x^{(\tau_3)}\oplus
x^{(\tau_4)}\oplus \cdots\oplus x^{(\tau_n)}$$ since it does not
matter where the parentheses lie.
\end{remark}

\subsection{Additive inverses}\label{subsection Golden Flowers}

\begin{defn}\label{defn41}
For any element $x=x_0.x_1x_2x_3\cdots\in{\mathcal R}$, let
$$\Psi(x)\triangleq(-1-x_0).(9-x_1)(9-x_2)(9-x_3)\cdots,$$
understood as  the ``additive inverse" of $x$. The absolute value of
$x$ to is defined to be
$$\| x\|\triangleq\max\{x,\Psi(x)\}.$$ Let $\mbox{sign}:{\mathcal
R}\rightarrow\{0,1\}$ be the signal map
\begin{equation*} \mbox{sign}(x)\triangleq
\begin{cases}
\ 0 \ \ \mbox{if}\  x\succeq0.000000\cdots, \\
\ 1 \ \ \mbox{if}\  x\preceq(-1).999999\cdots.
\end{cases}
\end{equation*}
\end{defn}

Some elementary properties on $\Psi(\cdot)$, $\|\cdot\|$ and
$\mbox{sign}(\cdot)$ are collected below without proofs. The
interested readers can easily provide the details without much
difficulty.

\begin{enumerate}

\item  $x\oplus\Psi(x)=(-1).999999\cdots$;
\item $\Psi(\Psi(x))=x$;
\item $x\preceq y\Leftrightarrow \Psi(x)\succeq\Psi(y)$;
\item $x\sim y\Leftrightarrow \Psi(x)\sim\Psi(y)$;

 \item $x\sim y\Rightarrow\|x\|\sim\|y\|$;
\item $x\sim\Psi(x)\Leftrightarrow x\sim0.000000\cdots$;
  \item   $\| x\|\succeq0.000000\cdots$;
  \item $\|\Psi(x)\|=\| x\|$;
\item
${\mbox{sign}}(x)+{\mbox{sign}}(\Psi(x))=1$;

\item $\Psi^{({\mbox{sign}}(x))}(x)=\|x\|$;
\item $\Psi^{({\mbox{sign}}(x))}(\|x\|)=x$;


\item   $\Psi(\sup W)=\inf \Psi(W)$;
\item
$\Psi(\overline{\mbox{LIMIT}}(\{(x^{(n)})\}_{n\in\mathbb{N}}))=\underline{\mbox{LIMIT}}(\{\Psi(x^{(n)})\}_{n\in\mathbb{N}})$.

\end{enumerate}

\begin{theorem}
Given $x,y,z\in{\mathcal R}$ with $x\oplus z\sim y\oplus z$, we have
$x\sim y$.
\end{theorem}

\begin{proof}
By Theorems \ref{keep order 33} and \ref{zzzzz},
\begin{align*}x&=x\oplus0.000000\cdots\\&\sim x\oplus z\oplus \Psi(z)\\&\sim y\oplus
z\oplus \Psi(z)\\&\sim y\oplus0.000000\cdots\\&=y.\end{align*}
\end{proof}

\begin{theorem}\label{abbbbbbbbbaaaaaaa}
For any two elements $x,y\in{\mathcal R}$, we have $\Psi(x\oplus
y)\sim\Psi(x)\oplus\Psi(y)$.
\end{theorem}

\begin{proof}
By Theorems \ref{abel sum}, \ref{keep order 33} and \ref{zzzzz},
\begin{align*}
\Psi(x\oplus y)&=\Psi(x\oplus
y)\oplus0.000000\cdots\oplus0.000000\cdots\\
&\sim \Psi(x\oplus y)\oplus(x\oplus\Psi(x))\oplus(y\oplus\Psi(y))\\
&\sim \Psi(x\oplus y)\oplus (x\oplus y)\oplus
(\Psi(x)\oplus\Psi(y))\\
&\sim\Psi(x)\oplus\Psi(y).
\end{align*}
\end{proof}

\section{Multiplicative operations}

\subsection{Multiplicative operations}

\begin{defn}[Multiplication]\label{defn51}
For any two elements $x,y\in{\mathcal R}$, let
\[x\otimes y\triangleq\Psi^{\big(\mbox{sign}(x)+\mbox{sign}(y)\big)}\Big(\mbox{LIMIT}\big(\{[\|x\|]_k\cdot[\|y\|]_k\}_{k\in\mathbb{N}}\big)\Big),\]
where $\Psi^{(k)}$ is the $k$-times composites of $\Psi$.
\end{defn}

\begin{example}\label{example52}
For any $x\in{\mathcal R}$,
$$x\otimes1.000000\cdots=\Psi^{(\mbox{sign}(x))}\big(\mbox{LIMIT}\big(\{[\|x\|]_k\}_{k\in\mathbb{N}}\big)\big)=\Psi^{({\mbox{sign}}(x))}(\|x\|)=x.$$
This means $({\mathcal R},\otimes)$ has a unit.
\end{example}

\begin{example}\label{example}
For any  $x\in{\mathcal R}$,
\begin{align*}x\otimes0.000000\cdots&=\Psi^{\big(\mbox{sign}(x)\big)}\Big(\mbox{LIMIT}\big(\{[\|x\|]_k\cdot0\}_{k\in\mathbb{N}}\big)\Big)
\\&=\Psi^{\big(\mbox{sign}(x)\big)}
(0.000000\cdots)\\&\sim0.000000\cdots,\\
x\otimes(-1).999999\cdots&=\Psi^{\big(\mbox{sign}(x)+1\big)}\Big(\mbox{LIMIT}\big(\{[\|x\|]_k\cdot0\}_{k\in\mathbb{N}}\big)\Big)
\\&=\Psi^{\big(\mbox{sign}(x)+1\big)}
(0.000000\cdots)\\&\sim0.000000\cdots.\end{align*}
\end{example}

\begin{example}
Given a calculator with sufficiently long digits, we could observe
that
\begin{align*}
1<1.4^2<1.99&<2<1.5^2\\
1.9<1.41^2<1.9999&<2<1.42^2\\
1.99<1.414^2<1.999999&<2<1.415^2\\
& \ \ \ \ \vdots\\
1.\underbrace{99\cdots99}_{n-1}<(a_0.a_1a_2a_3\cdots
a_n)^2<1.\underbrace{999\cdots999}_{2n}&<2<(a_0.a_1a_2a_3\cdots
a_n+\frac{1}{10^n})^2\\
&\ \ \ \ \vdots
\end{align*}
For thousands of years $a\triangleq a_0.a_1a_2a_3\cdots$ has been
understood as the positive square root of $2$, so what is the reason
behind? According to Definition \ref{defn51}, $a\otimes
a=1.999999\cdots$. Also from the above formulas, it is no hard to
observe (see also \cite{Courant,Flannery,Gowers0}) that there is no
element $z\in{\mathcal R}$ such that $z\otimes z=2.000000\cdots$. So
if we want to define the positive square root of $2$, except
$a_0.a_1a_2a_3\cdots$, which else could be?

\end{example}

\begin{remark}
Given $x,y\in\mathbb{Q}_F\subset\mathbb{Q}$ with
$x,y\succeq0.000000\cdots$, it is easy to verify that $x\otimes
y=x\cdot y$. Therefore from now on we can abuse the uses of
$\otimes$ and $\cdot$ if the summands lie in $\mathbb{Q}_F$ with
signs zero.
\end{remark}

\begin{theorem}\label{abel product} For any two elements $x,y\in{\mathcal R}$, we have
$x\otimes y=y\otimes x.$
\end{theorem}

\begin{proof}
\begin{align*}
x\otimes
y&=\Psi^{\big(\mbox{sign}(x)+\mbox{sign}(y)\big)}\Big(\mbox{LIMIT}\big(\{[\|x\|]_k\cdot[\|y\|]_k\}_{k\in\mathbb{N}}\big)\Big)\\
&=\Psi^{\big(\mbox{sign}(y)+\mbox{sign}(x)\big)}\Big(\mbox{LIMIT}\big(\{[\|y\|]_k\cdot[\|x\|]_k\}_{k\in\mathbb{N}}\big)\Big)=y\otimes
x.
\end{align*}
\end{proof}

\begin{theorem}\label{aaaaaaaa} For any two elements $x,y\in{\mathcal R}$, we have
$$\Psi(x\otimes y)=\Psi(x)\otimes y= x\otimes \Psi(y).$$
\end{theorem}

\begin{proof}
This is the twin theorem of Theorem \ref{abbbbbbbbbaaaaaaa}. By
Theorem \ref{abel product} it suffices to prove $\Psi(x\otimes
y)=\Psi(x)\otimes y$. To this aim we note
\begin{align*}
\Psi(x)\otimes y&=\Psi^{\big(\mbox{sign}(\Psi(x))+\mbox{sign}(y)\big)}\Big(\mbox{LIMIT}\big(\{[\|\Psi(x)\|]_k\cdot[\|y\|]_k\}_{k\in\mathbb{N}}\big)\Big)\\
&=\Psi^{\big(1+\mbox{sign}(x)+\mbox{sign}(y)\big)}\Big(\mbox{LIMIT}\big(\{[\|x\|]_k\cdot[\|y\|]_k\}_{k\in\mathbb{N}}\big)\Big)\\
&=\Psi\Bigg(\Psi^{\big(\mbox{sign}(x)+\mbox{sign}(y)\big)}\Big(\mbox{LIMIT}\big(\{[\|x\|]_k\cdot[\|y\|]_k\}_{k\in\mathbb{N}}\big)\Big)\Bigg)\\
&=\Psi(x\otimes y),
\end{align*}
here we have used the fact that $\Psi^{(2)}$ is the identity map.
This proves the theorem.

\end{proof}

\begin{theorem}\label{keep order 2} Given $x,y,z,w\in{\mathcal R}$ with $0.000000\preceq x\preceq z$ and $0.000000\preceq y\preceq w$, we have
$x\otimes y\preceq z\otimes w.$
\end{theorem}

\begin{proof}
This is the twin theorem of Theorem \ref{keep order}. $\forall
k\in\mathbb{N}$ we have $0\leq[x]_k\leq[z]_k$ and
$0\leq[y]_k\leq[w]_k$, which yields
$[x]_k\cdot[y]_k\leq[z]_k\cdot[w]_k$. By Lemma \ref{lemma order},
$\mbox{LIMIT}\big(\{[x]_k\cdot[y]_k\}_{k\in\mathbb{N}}\big)\preceq\mbox{LIMIT}\big(\{[z]_k\cdot
[w]_k\}_{k\in\mathbb{N}}\big)$. This proves the theorem.
\end{proof}

\begin{theorem}\label{keep order 333333} Given $x,y,z,w\in{\mathcal R}$ with $x\sim z$ and $y\sim w$, we have
$x\otimes y\sim z\otimes w.$
\end{theorem}

\begin{proof}
Obviously if two equivalent elements $x$ and $z$ have different
signs, then we must have
$\{x,z\}=\{0.000000\cdots,(-1).999999\cdots\}$. From Example
\ref{example}, $x\otimes y\sim0.000000\cdots\sim z\otimes w$. Thus
to prove this theorem, we may assume that
$\mbox{sign}(x)=\mbox{sign}(z)$, $\mbox{sign}(y)=\mbox{sign}(w)$,
which yields
$\mbox{sign}(x)+\mbox{sign}(y)=\mbox{sign}(z)+\mbox{sign}(w)$. By
Theorem \ref{aaaaaaaa},
\begin{align*}\Psi^{(\mbox{sign}(x)+\mbox{sign}(y))}(x\otimes
y)&=\Psi^{(\mbox{sign}(x))}(x)\otimes\Psi^{(\mbox{sign}(y))}(y)\ =\|x\|\otimes\|y\|,\\
\Psi^{(\mbox{sign}(z)+\mbox{sign}(w))}(z\otimes
w)&=\Psi^{(\mbox{sign}(z))}(z)\otimes\Psi^{(\mbox{sign}(w))}(w)=\|z\|\otimes\|w\|.\end{align*}
Consequently, to prove this theorem we may further assume below
$x,y,z,w\succeq0.000000\cdots$, by which we shall make use of
Theorem \ref{keep order 2}. Let $M\in\mathbb{N}$ be an upper bound
for $\{x,y,z,w\}$. For any $k\in\mathbb{N}$,
\begin{align*}
[x\otimes y]_k&\preceq x\otimes
y\preceq([x]_k+\frac{1}{10^k})\cdot([y]_k+\frac{1}{10^k})\\&\leq([z]_k+\frac{2}{10^k})
\cdot([w]_k+\frac{2}{10^k})\leq[z]_k\cdot[w]_k+\frac{5M}{10^k}\\&\preceq
(z\otimes w)\oplus\frac{5M}{10^k}\preceq[z\otimes
w]_k+\frac{5M+1}{10^k},
\end{align*}
here we have used Lemma \ref{equivalence test} for $x\sim z$ and
$y\sim w$. By symmetry we can also have $[z\otimes w]_k\leq
[x\otimes y]_k+\frac{5M+1}{10^k}.$ Finally by Lemma \ref{eqvi test
improved}, $x\otimes y\sim z\otimes w.$ This concludes the whole
proof.

\end{proof}

\begin{theorem}\label{associate product} For any three elements $x,y,z\in{\mathcal R}$, we have
$(x\otimes y)\otimes z\sim x\otimes(y\otimes z).$
\end{theorem}

\begin{proof}
By Theorem \ref{aaaaaaaa}, to prove this theorem we may assume that
$x,y,z\succeq0.000000\cdots$,  by which we shall also make use of
Theorem \ref{keep order 2}. Let $M\in\mathbb{N}$ be an upper bound
for $\{x,y,z\}$.  For any $k\in\mathbb{N}$, \begin{align*}(x\otimes
y)\otimes
z&\preceq([x]_k+\frac{1}{10^k})\cdot([y]_k+\frac{1}{10^k})\cdot([z]_k+\frac{1}{10^k})
\\&\leq[x]_k\cdot[y]_k\cdot[z]_k+\frac{4M^2}{10^k}\preceq
w\oplus\frac{4M^2}{10^k},
\end{align*}
where
$w\triangleq\mbox{LIMIT}\big(\{[x]_n\cdot[y]_n\cdot[z]_n\}_{n\in\mathbb{N}}\big)$.
 On the other hand,
$$[x]_n\cdot[y]_n\cdot[z]_n\preceq (x\otimes
y)\otimes[z]_n\preceq(x\otimes y)\otimes z.$$ Letting
$n\rightarrow\infty$ yields $w\preceq(x\otimes y)\otimes z$. With
these preparations in hand, for any $k\in\mathbb{N}$ we have
$[w]_k\leq[(x\otimes y)\otimes z]_k\leq
[w\oplus\frac{4M^2}{10^k}]_k=[w]_k+\frac{4M^2}{10^k}.$ Similarly, we
can also have $[w]_k\leq [(y\otimes z)\otimes x]_k\leq
[w]_k+\frac{4M^2}{10^k}$. Finally by Lemma \ref{eqvi test improved},
 $(x\otimes y)\otimes
z\sim(y\otimes z)\otimes x=x\otimes(y\otimes z)$. This concludes the
whole proof.
\end{proof}

\begin{theorem}\label{distribute law} For any three elements $x,y,z\in{\mathcal R}$, we have
$(x\oplus y)\otimes z\sim (x\otimes z)\oplus (y\otimes z).$
\end{theorem}

\begin{proof}
By Theorems \ref{abbbbbbbbbaaaaaaa} and \ref{aaaaaaaa}, to prove
this theorem we may assume $\mbox{sign}(z)=0$, and suppose this is
the case.

Case 1: Suppose $\mbox{sign}(x)=\mbox{sign}(y)$. By Theorems
\ref{abbbbbbbbbaaaaaaa} and \ref{aaaaaaaa} again, we may further
assume $\mbox{sign}(x)=\mbox{sign}(y)=0$. Let $M\in\mathbb{N}$ be an
upper bound for $\{x,y,z\}$. For any $k\in\mathbb{N}$,
\begin{align*}[(x\oplus y)\otimes z]_k&\preceq(x\oplus y)\otimes z\preceq([x\oplus
y]_k+\frac{1}{10^k})\cdot([z]_k+\frac{1}{10^k})\\
&\leq([x]_k+
[y]_k+\frac{3}{10^k})\cdot([z]_k+\frac{1}{10^k})\leq[x]_k\cdot[z]_k+[y]_k\cdot[z]_k+\frac{6M}{10^k}\\
&\preceq\big((x\otimes z)\oplus (y\otimes
z)\big)\oplus\frac{6M}{10^k}\preceq[(x\otimes z)\oplus (y\otimes
z)]_k+\frac{6M+1}{10^k}.\end{align*} On the other hand,
\begin{align*}[(x\otimes z)\oplus (y\otimes
z)]_k&\preceq(x\otimes z)\oplus (y\otimes
z)\\&\preceq([x]_k+\frac{1}{10^k})\cdot([z]_k+\frac{1}{10^k})+([y]_k+\frac{1}{10^k})\cdot([z]_k+\frac{1}{10^k})\\
&\leq([x]_k+[y]_k)\cdot[z]_k+\frac{6M}{10^k}\preceq((x\oplus
y)\otimes z)\oplus\frac{6M}{10^k}\\
&\preceq[(x\oplus y)\otimes z]_k+\frac{6M+1}{10^k}.\end{align*} Thus
by Lemma \ref{eqvi test improved}, $(x\oplus y)\otimes z\sim
(x\otimes z)\oplus (y\otimes z).$

Case 2: Suppose $x,y$ have different signs. Thus $x\oplus y$ has the
same sign as one of $x,y$, say for example, has the same sign as
$x$. Thus $\Psi(y)$ and $x\oplus y$ have the same sign. According to
the analysis in Case 1,
$$(\Psi(y)\oplus (x\oplus y))\otimes z\sim (\Psi(y)\otimes z)\oplus
((x\oplus y)\otimes z).$$ Adding $y\otimes z$ to the both sides of
the above formula yields
$$(y\otimes z)\oplus(x\otimes z)\sim(y\otimes
z)\oplus(\Psi(y)\otimes z)\oplus ((x\oplus y)\otimes z).$$ So we are
left to prove that $(y\otimes z)\oplus(\Psi(y)\otimes
z)\sim0.000000\cdots$. To this aim, we note
$$\Psi\big((y\otimes z)\oplus(\Psi(y)\otimes
z)\big)\sim(\Psi(y)\otimes z)\oplus(y\otimes z)=(y\otimes
z)\oplus(\Psi(y)\otimes z),$$ which naturally implies  $(y\otimes
z)\oplus(\Psi(y)\otimes z)\sim0.000000\cdots$.

This concludes the whole proof of the theorem.
\end{proof}

\subsection{Multiplicative inverses}

\begin{defn}[Inverse of multiplication]
For any $x\in{\mathcal R}$ with $\|x\|\succ0.000000\cdots$, let
\[x^{-1}\triangleq\Psi^{\big(\mbox{sign}(x)\big)}\Big(\mbox{LIMIT}\big(\{\frac{1}{[\|x\|]_k}\}_{k\in\mathbb{N}}\big)\Big).\]
For some small $k$, $[\|x\|]_k$ might be equal to zero. So it is
possible that $\frac{1}{[\|x\|]_k}$ is meaningless, or has the
meaning of positive infinity. But since we are mainly concerned with
the greatest lower bound for the set $\{\frac{1}{[\|x\|]_k}|\
k\in\mathbb{N}\}$, it does not matter whether the initial items of
the sequence  $\{\frac{1}{[\|x\|]_k}\}_{k\in\mathbb{N}}$ are equal
to positive infinity or not.
\end{defn}

\begin{theorem}\label{62}
For any $x\in{\mathcal R}$ with $\|x\|\succ0.000000\cdots$, we have
$\Psi(x)^{-1}=\Psi(x^{-1})$.
\end{theorem}

\begin{proof}
Since $\|\Psi(x)\|=\|x\|\succ0.000000\cdots$, by definition we have
\begin{align*}\Psi(x)^{-1}&=\Psi^{\big(\mbox{sign}(\Psi(x))\big)}\Big(\mbox{LIMIT}\big(\{\frac{1}{[\|\Psi(x)\|]_k}\}_{k\in\mathbb{N}}\big)\Big)
\\&=\Psi^{\big(1+\mbox{sign}(x)\big)}\Big(\mbox{LIMIT}\big(\{\frac{1}{[\|x\|]_k}\}_{k\in\mathbb{N}}\big)\Big)=\Psi(x^{-1}).\end{align*}
\end{proof}

\begin{theorem}\label{muuuuuuuuu}
For any $x\in{\mathcal R}$ with $\|x\|\succ0.000000\cdots$, we have
$x\otimes x^{-1}\sim1.000000\cdots$.
\end{theorem}

\begin{proof}
According to Theorem \ref{62}, to prove this theorem we may assume
$x\succ0.000000\cdots$, and suppose this is the case. Since
$\|x\|\succ0.000000\cdots$, one can find a sufficiently large
$m\in\mathbb{N}$ such that $\frac{1}{10^m}\preceq x$. For any
naturals number $k\geq n\geq m$, $[x]_n\leq[x]_k\leq[x]_n+10^{-n}$.
Hence taking reciprocals in each part we get
$\frac{1}{[x]_n+10^{-n}}\leq\frac{1}{[x]_k}\leq\frac{1}{[x]_n}$.
Letting $k\rightarrow\infty$ yields $\frac{1}{[x]_n+10^{-n}}\preceq
x^{-1}\preceq\frac{1}{[x]_n}$. Combining this with Lemma
\ref{additive estimate one element KEY} we have
$$[x]_n\cdot\frac{1}{[x]_n+\frac{1}{10^n}}\preceq x\otimes x^{-1}\preceq([x]_n+\frac{1}{10^n})\cdot\frac{1}{[x]_n}.$$
Since $\frac{1}{10^m}\preceq x$, we have
$\frac{1}{10^m}\leq[x]_m\leq[x]_n$, and consequently,
$[x]_n\cdot\frac{1}{[x]_n+\frac{1}{10^n}}\geq1-10^{m-n}$,
$([x]_n+\frac{1}{10^n})\cdot\frac{1}{[x]_n}\leq1+10^{m-n}.$ Finally
for any $s\in\mathbb{N}$, let us take $n=m+s$. Then
$1-\frac{1}{10^s}\preceq x\otimes x^{-1}\preceq1+\frac{1}{10^s}$,
which yields $1-\frac{1}{10^s}\leq [x\otimes
x^{-1}]_s\leq1+\frac{1}{10^s}$. Note also
$1-\frac{1}{10^s}\leq[1.000000\cdots]_s\leq1+\frac{1}{10^s}$. Thus a
standard application of Lemma \ref{equivalence test} gives the
desired result. We are done.

\end{proof}

\section{Construction of real numbers}

\subsection{Motivation} In mathematics, since ${\mathcal R}$ is a
set, there is no doubt that $0.999999\cdots$ and $1.000000\cdots$
are two different elements. But in the real world, if we still
imagine $0.999999\cdots$ and $1.000000\cdots$ as two different
objects, then it will cause some trouble. For example, how can we
define and understand
$$\frac{0.999999\cdots+1.000000\cdots}{2}?$$
Since there is no element $z\in{\mathcal R}$ such that
$0.999999\cdots\prec z\prec1.000000\cdots$, we may regard ${\mathcal
R}$ having ``holes" in many places. Have you ever seen a bunch of
sunlight having lots of holes? It is too terrible!

Noting the sequence $\{[0.999999\cdots]_k\}_{k=1}^{\infty}$ is
getting ``closer and closer" to $1.000000\cdots$ as $k$ approaches
to infinity,  the best way we think to understand the relation
between $0.999999\cdots$ and $1.000000\cdots$ both in mathematics
and in the real world, is to view them as the same object, which
leads to the invention of the real ``real numbers" below.

\subsection{Methodology} Define the set of real numbers $\mathbb{R}$
to be the set of equivalent classes ${\mathcal R}/\sim$. For those
readers who are not too familiar with the algebraic representation
${\mathcal R}/\sim$, the definition of $\mathbb{R}$ is not too
abstract at all: if the equivalent class has exactly two elements,
you may simply discard one element, say for example discard those
elements in ${\mathcal R}_9$ as we have done in Section
\ref{section1}, and leave the other one remained.

For any two equivalent classes $[x],[y]$, we introduce an induced
operation $\widehat{\oplus}$ on $\mathbb{R}$ by
\begin{align*}
[x]\widehat{\oplus}[y]&\triangleq[x\oplus y].
\end{align*}
By Theorems \ref{keep order 33}, $\widehat{\oplus}$ is well-defined.
By Theorem \ref{abel sum}, $\widehat{\oplus}$ is commutative. By
Theorem \ref{zzzzz}, $\widehat{\oplus}$ is associative. From Example
\ref{example32}, $(\mathbb{R},\widehat{\oplus})$ has a unit
$[0.000000\cdots]$. Since
$x\oplus\Psi(x)=(-1).999999\cdots\sim0.000000\cdots$, every element
of $(\mathbb{R},\widehat{\oplus})$ has an inverse. In the language
of algebra,   $(\mathbb{R},\widehat{\oplus})$ is an Abelian group.

For any two equivalent classes $[x],[y]$, we introduce an induced
operation $\widehat{\otimes}$ on $\mathbb{R}$ by
\begin{align*}[x]\widehat{\otimes}[y]&\triangleq[x\otimes y].\end{align*}
By Theorem \ref{keep order 333333}, $\widehat{\otimes}$ is
well-defined. By Theorem \ref{abel product}, $\widehat{\otimes}$ is
commutative. By Theorem \ref{associate product}, $\widehat{\otimes}$
is associative. From Example \ref{example52},
$(\mathbb{R},\widehat{\otimes})$ has a unit $[1.000000\cdots]$. By
Theorem \ref{muuuuuuuuu}, every element of
$(\mathbb{R}^{\ast},\widehat{\otimes})$ has an inverse, where
$\mathbb{R}^{\ast}\triangleq\mathbb{R}\backslash[0.000000\cdots]$.
By Theorem \ref{distribute law}, $\widehat{\otimes}$ is distribute
over $\widehat{\oplus}$. Hence we have derived a number system
$(\mathbb{R},\widehat{\oplus},\widehat{\otimes})$, which in the
language of algebra is a field.

We can also introduce an order $\widehat{\preceq}$ on $\mathbb{R}$
from the total order $\preceq$ on ${\mathcal R}$, that is, we say
$[x]\widehat{\preceq}[y]$ if $y\oplus
\Psi(x)\succeq(-1).999999\cdots$. The interested readers can verify
that $\precsim$ is well-defined and is a total order. Finally, we
state the least upper  and greatest lower bounds properties for
$(\mathbb{R},\widehat{\preceq})$ below without proofs, which is a
simple corollary of Theorem \ref{main result}.

\begin{theorem}
Every nonempty bounded above (below) subset of
$(\mathbb{R},\widehat{\preceq})$ has a least upper (greatest lower)
bound.
\end{theorem}

As experienced readers should know, from the least upper bound
property, after introducing Weierstrass's $\epsilon-N$ definition of
convergence of sequences, we can deduce many other theorems such as
the monotone convergence theorem, Cauchy's convergence principle,
the Bolzano-Weierstrass theorem, the Heine-Borel theorem and so on
in a few pages.

\section{From decimals to Dedekind cuts}\label{sectiondedekind}

\subsection{A natural bijection}
\begin{defn}
A \textsf{Dedekind cut} is a pair of nonempty subsets $A,B$ of
$\mathbb{Q}$, denoted by $(A|B)$, such that
\begin{itemize}
\item $A\cap B=\emptyset$, $A\cup B=\mathbb{Q}$,
\item $a\in A, b\in B\Rightarrow a<b$,
\item $A$ contains no greatest element.
\end{itemize}
The set of all Dedekind cuts is denoted by $\mathbf{DR}$.
\end{defn}

Given a Dedekind cut $(A|B)$, we obviously have $\sup A\preceq\inf
B$. To go even further, we can have the following theorem.

\begin{theorem}\label{dedekind}
For any Dedekind cut $(A|B)$, we have $\sup A\sim\inf B$.
\end{theorem}

\begin{proof}
Fix arbitrarily $k\in\mathbb{N}$, $a\in A$ and $b\in B$. Step 1: Let
$c^{(1)}\triangleq\frac{a+b}{2}$. If $c^{(1)}\in A$, then we let
$a^{(1)}\triangleq c^{(1)}$, $b^{(1)}\triangleq b$; else if
$c^{(1)}\in B$, then we let $a^{(1)}\triangleq a$,
$b^{(1)}\triangleq c^{(1)}$. Obviously
$b^{(1)}-a^{(1)}=\frac{b-a}{2}$. Step 2: Let
$c^{(2)}\triangleq\frac{a^{(1)}+b^{(1)}}{2}$. If $c^{(2)}\in A$,
then we let $a^{(2)}\triangleq c^{(2)}$, $b^{(2)}\triangleq
b^{(1)}$; else if $c^{(2)}\in B$, then we let $a^{(2)}\triangleq
a^{(1)}$, $b^{(2)}\triangleq c^{(2)}$. Obviously
$b^{(2)}-a^{(2)}=\frac{b^{(1)}-a^{(1)}}{2}=\frac{b-a}{4}$. We can
repeat this procedure to the $4k$-th step to get $a^{(4k)}\in A$ and
$b^{(4k)}\in B$ so that $b^{(4k)}-a^{(4k)}=\frac{b-a}{16^k}$. Note
$a^{(4k)}\preceq \sup A\preceq \inf B\preceq b^{(4k)}$, thus
\[0\leq[\inf B]_k-[\sup A]_k\leq [b^{(4k)}]_k-[a^{(4k)}]_k\leq  b^{(4k)}-(a^{(4k)}-\frac{1}{10^k})
\leq\frac{\lceil b-a\rceil}{16^k}+\frac{1}{10^k}\leq\frac{2\lceil
b-a\rceil}{10^k},\] here  we have used Lemma \ref{additive estimate
one element KEY} at the third inequality. By Lemma \ref{eqvi test
improved}, $\sup A\sim\inf B$.
\end{proof}

 Now we can derive a
map $\tau$ from $\mathbf{DR}$ to $\mathcal{R}/\sim$ by sending
$(A|B)$ to $[\sup A]=[\inf B]$. It would be great if $\tau$ is
bijective. In the following we shall prove this is indeed the case.

\begin{lemma}\label{Interesting} We have
${\mathcal R}_9\cap \mathbb{Q}=\emptyset$.
\end{lemma}

\begin{proof}
Suppose the contrary one can find an element $\alpha\in {\mathcal
R}_9\cap \mathbb{Q}$. Let us choose the corresponding element
$\beta\in\mathbb{Q}_F$ so that $\alpha\prec\beta$,
$\alpha\sim\beta$. Obviously as elements of $\mathbb{Q}$,
$\alpha<\beta$. Now choose a sufficiently large $k\in \mathbb{N}$ so
that $\beta-\alpha\geq\frac{3}{10^k}$, then by Lemma \ref{additive
estimate one element KEY},
$$[\beta]_k-[\alpha]_k\geq (\beta-\frac{1}{10^k})-\alpha\geq\frac{2}{10^k}.$$
According to Lemma \ref{equivalence test}, this is impossible. We
are done.

\end{proof}

\begin{theorem}\label{dedekindbijection}
The map  $\tau$ from $\mathbf{DR}$ to
$\mathcal{R}/\sim$ is bijective.
\end{theorem}

\begin{proof}
We first prove that $\tau$ is surjective. Fix arbitrarily
$x\in\mathcal{R}$. Let
\begin{align*}
A&\triangleq\{q\in\mathbb{Q}|\ q\prec x\},\\
B&\triangleq\{q\in\mathbb{Q}|\ q\succeq x\}.
\end{align*}
Obviously to verify that $(A|B)$ is a Dedekind cut, it suffices to
prove that $A$ contains no greatest element. Suppose the contrary,
then $\sup A\in A$, which by the definitions of $A$ and $B$ leads to
$\sup A\prec x\preceq \inf B.$ By Theorem \ref{dedekind}, $\sup
A\sim\inf B$. By Lemma \ref{R9Qf}, $\sup A\in {\mathcal R}_9$, which
contradicts to Lemma \ref{Interesting}. Hence $(A|B)$ must be a
Dedekind cut. Note $\sup A\preceq x\preceq \inf B$, by Theorem
\ref{dedekind} again, we have $\tau((A|B))=[x]$. Thus $\tau$ is
surjective.

Next we prove that $\tau$ is injective. Given two different Dedekind
cuts $(A|B)$,  $(C|D)$, without loss of generality we may assume
there exists an element $\alpha\in C\backslash A$, which by the
definition of Dedekind cut yields $\sup A\preceq\alpha\prec\sup C$.
Our aim is to prove that $\sup A$ and $\sup C$ are not equivalent.
Suppose this is the contrary, then $\sup A\sim\sup C$. By Lemma
\ref{R9Qf}, $\alpha=\sup A\in{\mathcal R}_9$, which contradicts to
Lemma \ref{Interesting}. Thus $\tau$ is injective.

This concludes the whole proof of the theorem.
\end{proof}

\subsection{Isomorphisms between operations}

\begin{defn}A \textsf{generalized Dedekind cut} is a pair of nonempty
subsets $A,B$ of $\mathbb{Q}$, denoted by $[A\Rsh B]$, such that
\begin{itemize}
\item $A\cap B=\emptyset$, $A\cup B=\mathbb{Q}$,
\item $a\in A, b\in B\Rightarrow a<b$.
\end{itemize}
\end{defn}

\begin{remark}
To form a Dedekind cut $(\widetilde{A}|\widetilde{B})$ from a
generalized Dedekind cut $[A\Rsh B]$, we can simply move the
greatest element of $A$ if it exists, from $A$ to $B$. In this case
it is no hard to prove that $\sup A\sim\inf B\sim\sup
\widetilde{A}\sim\inf \widetilde{B}$, so we can abuse the uses of
$[A\Rsh B]$ and $(\widetilde{A}|\widetilde{B})$. For example, let
$\widetilde{\Psi}:\mathbf{DR}\rightarrow\mathbf{DR}$ be the additive
inverse map defined by sending $(A|B)$ to $[-B\Rsh -A]$.
\end{remark}

\begin{defn} Let
$\Theta:\mathbf{DR}\rightarrow\{0,1\}$ be the signal map
\[\Theta((A|B))\triangleq \max\limits_{x\in B}\mbox{sign}(x).\] The
absolute value of a Dedekind cut $(A|B)$ is defined by
\[\|(A|B)\|\triangleq\widetilde{\Psi}^{\big(\Theta((A|B))\big)}((A|B)).\]
\end{defn}

\begin{defn}

For any two Dedekind cuts $(A|B)$, $(C|D)$, define the addition of
both cuts by
\[(A|B)\boxplus(C|D)=[\mathbb{Q}\backslash(B+D)\Rsh B+D].\] If the Dedekind cuts $(A|B)$,
$(C|D)$ are of signs zero, define the multiplication of them by
\[(A|B)\boxtimes(C|D)=[\mathbb{Q}\backslash(B\cdot D)\Rsh B\cdot D].\]
Generally for any two Dedekind cuts $(A|B)$, $(C|D)$, define  their
multiplication by
\[(A|B)\boxtimes(C|D)=\widetilde{\Psi}^{\big(\Theta((A|B))+\Theta((C|D))\big)}\Big(\|(A|B)\|\boxtimes\|(C|D)\|\Big).\]

\end{defn}

\begin{theorem}\label{dedekindadditive}
For any two Dedekind cuts $(A|B)$, $(C|D)$, we have
\[\tau\big((A|B)\boxplus(C|D)\big)=\tau\big((A|B)\big)\widehat{\oplus}\ \tau\big((C|D)\big).\]
\end{theorem}

\begin{proof}
To prove this theorem, it suffices to prove that $\inf (B+D)\sim
\inf B \oplus \inf D$. To this aim we  note $\sup(A+C)\preceq \sup
A\oplus \sup C\preceq \inf B \oplus \inf D\preceq \inf (B+D).$ Hence
we need only to prove that $\sup(A+C)\sim  \inf (B+D)$.

Fix arbitrarily $k\in\mathbb{N}$. We can repeat the procedures in
the proof of Theorem \ref{dedekind} from any common starting points
$a\in A\cap C$ and $b\in B\cap D$ to get $a^{(4k)}\in A$,
$b^{(4k)}\in B$, $c^{(4k)}\in C$ and $d^{(4k)}\in D$ so that
$b^{(4k)}-a^{(4k)}=d^{(4k)}-c^{(4k)}=\frac{b-a}{16^k}$. Thus
\[a^{(4k)}+ c^{(4k)}\preceq \sup (A+C)\preceq \inf (B+ D)\preceq
b^{(4k)}+ d^{(4k)},\] and consequently,
 \[0\leq[\inf (B+ D)]_k-[\sup (A+
C)]_k\leq[ b^{(4k)}+ d^{(4k)}]_k-[a^{(4k)}+
c^{(4k)}]_k\leq\frac{2\lceil b-a\rceil}{16^k}+\frac{1}{10^k}.\] By
Lemma \ref{eqvi test improved}, we are done.
\end{proof}

\begin{theorem}\label{iso}
For any two Dedekind cuts $(A|B)$, $(C|D)$, we have
\[\tau\big((A|B)\boxtimes(C|D)\big)=\tau\big((A|B)\big)\widehat{\otimes}\ \tau\big((C|D)\big).\]

\end{theorem}

\begin{proof}
We shall only prove this theorem for the special case that the
Dedekind cuts $(A|B)$, $(C|D)$ are of signs zero, and leave the
general case to the interested readers. Suppose the Dedekind cuts
$(A|B)$, $(C|D)$ are of signs zero. To our aim, it suffices to prove
that $\inf (B\cdot D)\sim \inf B \otimes \inf D$. If either $0\in B$
or $0\in D$, then it is easy to observe that $\inf (B\cdot D)=\inf B
\otimes \inf D=0.000000\cdots.$ Hence we may further assume that
$0\in A\cap C$. Let us choose a fixed natural number $z\in B\cap D$.

Fix arbitrarily $k\in\mathbb{N}$. We can repeat the procedures in
the proof of Theorem \ref{dedekind} from the starting points $0\in
A\cap C$ and $z\in B\cap D$ to get $a^{(4k)}\in A$, $b^{(4k)}\in B$,
$c^{(4k)}\in C$ and $d^{(4k)}\in D$ so that
$b^{(4k)}-a^{(4k)}=d^{(4k)}-c^{(4k)}=\frac{z-0}{16^k}=\frac{z}{16^k}$.
Thus \[0.000000\cdots\preceq a^{(4k)}\cdot c^{(4k)}\preceq \sup A
\otimes \sup C\preceq \inf B \otimes \inf D\preceq \inf (B\cdot
D)\preceq b^{(4k)}\cdot d^{(4k)},\] and consequently, \[0\leq[\inf
(B\cdot D)]_k-[\inf B \otimes \inf D]_k\leq[ b^{(4k)}\cdot
d^{(4k)}]_k-[a^{(4k)}\cdot
c^{(4k)}]_k\leq\frac{2z^2}{16^k}+\frac{z^2}{256^k}+\frac{1}{10^k}.\]
By Lemma \ref{eqvi test improved}, we are done.
\end{proof}

Given two Dedekind cuts $(A|B)$, $(C|D)$, we say $(A|B)$ is less
than or equal to $(C|D)$, denoted by $(A|B)\lesssim(C|D)$, if $\sup
A\preceq \sup C$. It is a piece of cake to prove that
$$(A|B)\lesssim(C|D)\Leftrightarrow \tau((A|B))\ \widehat{\preceq}\ \tau((C|D)).$$
Thus with Theorems \ref{dedekindbijection}, \ref{dedekindadditive}
and \ref{iso},
$(\mathbb{R},\widehat{\oplus},\widehat{\otimes},\widehat{\preceq})$
is isomorphic to $(\mathbf{DR},\boxplus,\boxtimes,\lesssim)$.

\section{From decimals to Cauchy sequences}\label{Cauchy-Decimals}

\subsection{A natural surjection}
\begin{defn}
A  sequence of rational numbers $\{x^{(n)}\}_{n\in\mathbb{N}}$ is
called \textsf{Cauchy} if for every  rational  $\epsilon>0$, there
exists an $N\in\mathbb{N}$ such that  $\forall m,n \geq N$,
$|x^{(m)}-x^{(n)}|<\epsilon$. The set of all Cauchy sequences of
rational numbers is denoted by $\mathbf{CR}$.  Two Cauchy sequences
of rational numbers $\{x^{(n)}\}_{n\in\mathbb{N}}$,
$\{y^{(n)}\}_{n\in\mathbb{N}}$ are said to be equivalent, denoted by
$\{x^{(n)}\}_{n\in\mathbb{N}}\approx\{y^{(n)}\}_{n\in\mathbb{N}}$,
if for every  rational $\epsilon>0$, there exists an
$N\in\mathbb{N}$ such that $\forall n \geq N$,
$|x^{(n)}-y^{(n)}|<\epsilon$.
\end{defn}

It is easy to prove $\approx$ is an equivalence relation, thus we
can get a quotient space $\mathbf{CR}/\approx$, denoted by
$\widehat{\mathbf{CR}}$ for consistency. Given a Cauchy sequence of
rational numbers $\{x^{(n)}\}_{n\in\mathbb{N}}$, we obviously have
$\underline{\mbox{LIMIT}}(\{x^{(n)}\}_{n\in\mathbb{N}})\preceq\overline{\mbox{LIMIT}}(\{x^{(n)}\}_{n\in\mathbb{N}})$.
To go even further, we can have the following theorem.

\begin{theorem}\label{cauchy}
$\underline{\mbox{LIMIT}}(\{x^{(n)}\}_{n\in\mathbb{N}})\sim\overline{\mbox{LIMIT}}(\{x^{(n)}\}_{n\in\mathbb{N}})$
holds for any Cauchy   sequence of rational numbers
$\{x^{(n)}\}_{n\in\mathbb{N}}$.
\end{theorem}

\begin{proof} Fix arbitrarily $k\in\mathbb{N}$, there exists an $N\in\mathbb{N}$ depending on $k$ such that $\forall m,n \geq
N$, $|x^{(m)}-x^{(n)}|<\frac{1}{10^k}$. Thus $\forall n\geq N$,
$x^{(N)}-\frac{1}{10^k}\leq x^{(n)}\leq x^{(N)}+\frac{1}{10^k}$.
Letting $n\rightarrow\infty$ gives
\[x^{(N)}-\frac{1}{10^k}\preceq
\underline{\mbox{LIMIT}}(\{x^{(n)}\}_{n\in\mathbb{N}})\preceq
\overline{\mbox{LIMIT}}(\{x^{(n)}\}_{n\in\mathbb{N}})\preceq
x^{(N)}+\frac{1}{10^k},\] and consequently
\[0
\leq[\overline{\mbox{LIMIT}}(\{x^{(n)}\}_{n\in\mathbb{N}})]_k-[\underline{\mbox{LIMIT}}(\{x^{(n)}\}_{n\in\mathbb{N}})]_k\leq
[x^{(N)}+\frac{1}{10^k}]_k-[x^{(N)}-\frac{1}{10^k}]_k=\frac{2}{10^k}.\]
By Lemma \ref{eqvi test improved},
$\underline{\mbox{LIMIT}}(\{x^{(n)}\}_{n\in\mathbb{N}})\sim\overline{\mbox{LIMIT}}(\{x^{(n)}\}_{n\in\mathbb{N}})$.
This finishes the proof.

\end{proof}
Now we can derive a map $\kappa$ from $\mathbf{CR}$ to
$\mathcal{R}/\sim$ by sending Cauchy sequence
$\{x^{(n)}\}_{n\in\mathbb{N}}$ to
$[\underline{\mbox{LIMIT}}(\{x^{(n)}\}_{n\in\mathbb{N}})]=[\overline{\mbox{LIMIT}}(\{x^{(n)}\}_{n\in\mathbb{N}})]$.
Obviously, $\kappa$ is surjective since for any $x\in{\mathcal R}$,
we have a typical Cauchy sequence $\{[x]_k\}_{k\in\mathbb{N}}$ such
that $x=\mbox{LIMIT}\big(\{[x]_k\}_{k\in\mathbb{N}}\big)$.

\begin{theorem}\label{cauchy72}
For any two Cauchy sequences of rational numbers
$\{x^{(n)}\}_{n\in\mathbb{N}}$, $\{y^{(n)}\}_{n\in\mathbb{N}}$,
$\{x^{(n)}\}_{n\in\mathbb{N}}\approx\{y^{(n)}\}_{n\in\mathbb{N}}\Leftrightarrow
\kappa(\{x^{(n)}\}_{n\in\mathbb{N}})=\kappa(\{y^{(n)}\}_{n\in\mathbb{N}})$.
\end{theorem}

\begin{proof} Suppose
$\{x^{(n)}\}_{n\in\mathbb{N}}\approx\{y^{(n)}\}_{n\in\mathbb{N}}$.
Fix arbitrarily $k\in\mathbb{N}$, there exists an $N\in\mathbb{N}$
such that  $\forall n \geq N$, $|x^{(n)}-y^{(n)}|<\frac{1}{10^k}$.
Equivalently, $\forall n \geq N$ we have $x^{(n)}\leq
y^{(n)}+\frac{1}{10^k}$ and $y^{(n)}\leq x^{(n)}+\frac{1}{10^k}$.
Letting $n\rightarrow\infty$ yields
$\overline{\mbox{LIMIT}}(\{x^{(n)}\}_{n\in\mathbb{N}})\preceq\overline{\mbox{LIMIT}}(\{y^{(n)}\}_{n\in\mathbb{N}})\oplus\frac{1}{10^k}$
and
$\overline{\mbox{LIMIT}}(\{y^{(n)}\}_{n\in\mathbb{N}})\preceq\overline{\mbox{LIMIT}}(\{x^{(n)}\}_{n\in\mathbb{N}})\oplus\frac{1}{10^k}$,
which further implies
\[\big|[\overline{\mbox{LIMIT}}(\{x^{(n)}\}_{n\in\mathbb{N}})]_k-[\overline{\mbox{LIMIT}}(\{y^{(n)}\}_{n\in\mathbb{N}})]_k\big|\leq\frac{1}{10^k}.\]
By Lemma \ref{eqvi test improved},
$\overline{\mbox{LIMIT}}(\{x^{(n)}\}_{n\in\mathbb{N}})\sim\overline{\mbox{LIMIT}}(\{y^{(n)}\}_{n\in\mathbb{N}})$.
This proves the necessary part.

Suppose
$\kappa(\{x^{(n)}\}_{n\in\mathbb{N}})=\kappa(\{y^{(n)}\}_{n\in\mathbb{N}})$.
We argue by contradiction the sufficient part and suppose
$\{x^{(n)}\}_{n\in\mathbb{N}}$, $\{y^{(n)}\}_{n\in\mathbb{N}}$ are
not equivalent. Then there exist a positive rational number
$\epsilon_0$ and a sequence of natural numbers $n_1<n_2<n_3<\cdots$
such that $\forall i\in \mathbb{N}$,
$|x^{(n_i)}-y^{(n_i)}|\geq\epsilon_0$. It may happen either
$x^{(n_i)}\geq y^{(n_i)}+\epsilon_0$ or $y^{(n_i)}\geq
x^{(n_i)}+\epsilon_0$. Without loss of generality, we may assume
there are infinitely many terms $x^{(n_i)}\geq
y^{(n_i)}+\epsilon_0$. Now  we first choose a $k\in\mathbb{N}$ such
that $\frac{2}{10^k}\leq\epsilon_0$, then choose a subsequence
$\{n_{i_j}\}_{j\in\mathbb{N}}$ from $\{n_{i}\}_{i\in\mathbb{N}}$
such that $x^{(n_{i_j})}\geq y^{(n_{i_j})}+\frac{2}{10^k}$. Thus
\[\overline{\mbox{LIMIT}}(\{x^{(n_{i_j})}\}_{j\in\mathbb{N}})\succeq\overline{\mbox{LIMIT}}(\{y^{(n_{i_j})}\}_{j\in\mathbb{N}})\oplus\frac{2}{10^k},\]
which yields
$[\overline{\mbox{LIMIT}}(\{x^{(n_{i_j})}\}_{j\in\mathbb{N}})]_k\geq[\overline{\mbox{LIMIT}}(\{y^{(n_{i_j})}\}_{j\in\mathbb{N}})]_k+\frac{2}{10^k}$,
and consequently by Lemma \ref{eqvi test improved},
$\overline{\mbox{LIMIT}}(\{x^{(n_{i_j})}\}_{j\in\mathbb{N}})$ and
$\overline{\mbox{LIMIT}}(\{y^{(n_{i_j})}\}_{j\in\mathbb{N}})$ are
not equivalent in $\mathcal R$. But this is impossible since by
Lemma \ref{lemma order}, Theorem \ref{cauchy} and our assumptions,
\begin{align*}
\underline{\mbox{LIMIT}}(\{x^{(n)}\}_{n\in\mathbb{N}})&\preceq\underline{\mbox{LIMIT}}(\{x^{(n_{i_j})}\}_{j\in\mathbb{N}})\preceq\overline{\mbox{LIMIT}}(\{x^{(n_{i_j})}\}_{j\in\mathbb{N}})\preceq\overline{\mbox{LIMIT}}(\{x^{(n)}\}_{n\in\mathbb{N}}),\\
\underline{\mbox{LIMIT}}(\{y^{(n)}\}_{n\in\mathbb{N}})&\preceq\underline{\mbox{LIMIT}}(\{y^{(n_{i_j})}\}_{j\in\mathbb{N}})\preceq\overline{\mbox{LIMIT}}(\{y^{(n_{i_j})}\}_{j\in\mathbb{N}})\preceq\overline{\mbox{LIMIT}}(\{y^{(n)}\}_{n\in\mathbb{N}}),\\
\underline{\mbox{LIMIT}}(\{x^{(n)}\}_{n\in\mathbb{N}})&\sim\overline{\mbox{LIMIT}}(\{x^{(n)}\}_{n\in\mathbb{N}})\sim\underline{\mbox{LIMIT}}(\{y^{(n)}\}_{n\in\mathbb{N}})\sim\overline{\mbox{LIMIT}}(\{y^{(n)}\}_{n\in\mathbb{N}}),
\end{align*}
which naturally implies that
$\overline{\mbox{LIMIT}}(\{x^{(n_{i_j})}\}_{j\in\mathbb{N}})\sim\overline{\mbox{LIMIT}}(\{y^{(n_{i_j})}\}_{j\in\mathbb{N}})$.
This proves the sufficient part, also concludes the whole proof of
the theorem.
\end{proof}

\subsection{Homomorphisms between operations}

\begin{defn}

Given two Cauchy sequences of rational numbers
$\{x^{(n)}\}_{n\in\mathbb{N}}$, $\{y^{(n)}\}_{n\in\mathbb{N}}$, it
is easy to verify that $\{x^{(n)}+ y^{(n)}\}_{n\in\mathbb{N}}$ and
$\{x^{(n)}\cdot y^{(n)}\}_{n\in\mathbb{N}}$ are also Cauchy
sequences, so we can define the addition and multiplication
respectively by
\begin{align*}\{x^{(n)}\}_{n\in\mathbb{N}}\widehat{\boxplus}\ \{y^{(n)}\}_{n\in\mathbb{N}}&=\{x^{(n)}+y^{(n)}\}_{n\in\mathbb{N}},\\
\{x^{(n)}\}_{n\in\mathbb{N}}\widehat{\boxtimes}\
\{y^{(n)}\}_{n\in\mathbb{N}}&=\{x^{(n)}\cdot
y^{(n)}\}_{n\in\mathbb{N}}.
\end{align*}
 \end{defn}

\begin{theorem}\label{cauchyadditive}
For any two Cauchy sequences of rational numbers
$\{x^{(n)}\}_{n\in\mathbb{N}}$, $\{y^{(n)}\}_{n\in\mathbb{N}}$,
\[\kappa\big(\{x^{(n)}\}_{n\in\mathbb{N}}\widehat{\boxplus}\
\{y^{(n)}\}_{n\in\mathbb{N}}\big)=
\kappa\big(\{x^{(n)}\}_{n\in\mathbb{N}}\big)\widehat{\oplus}\
\kappa\big(\{y^{(n)}\}_{n\in\mathbb{N}}\big).\]
\end{theorem}

\begin{proof}

 Fix arbitrarily $k\in\mathbb{N}$, there exists a common starting index $N\in\mathbb{N}$ such that $\forall m,n \geq
N$, $|x^{(m)}-x^{(n)}|<\frac{1}{10^k}$,
$|y^{(m)}-y^{(n)}|<\frac{1}{10^k}$. Thus $\forall n\geq N$,
\begin{align*}
x^{(N)}-\frac{1}{10^k}&\leq x^{(n)}\leq x^{(N)}+\frac{1}{10^k},\\
y^{(N)}-\frac{1}{10^k}&\leq y^{(n)}\leq y^{(N)}+\frac{1}{10^k},\\
x^{(N)}+y^{(N)}-\frac{2}{10^k}&\leq x^{(n)}+ y^{(n)}\leq
x^{(N)}+y^{(N)}+\frac{2}{10^k}.
\end{align*}
Letting $n\rightarrow\infty$ gives
\begin{align}
x^{(N)}-\frac{1}{10^k}&\preceq \underline{\mbox{LIMIT}}(\{x^{(n)}\}_{n\in\mathbb{N}})\preceq x^{(N)}+\frac{1}{10^k},\label{71}\\
y^{(N)}-\frac{1}{10^k}&\preceq \underline{\mbox{LIMIT}}(\{y^{(n)}\}_{n\in\mathbb{N}})\preceq y^{(N)}+\frac{1}{10^k},\label{72}\\
x^{(N)}+y^{(N)}-\frac{2}{10^k}&\preceq
\underline{\mbox{LIMIT}}(\{x^{(n)}+y^{(n)}\}_{n\in\mathbb{N}})\preceq
x^{(N)}+y^{(N)}+\frac{2}{10^k}.\label{73}
\end{align}
Applying Theorem \ref{keep order} to (\ref{71}) and (\ref{72})
yields
\begin{align}\label{74}x^{(N)}+y^{(N)}-\frac{2}{10^k}\preceq \underline{\mbox{LIMIT}}(\{x^{(n)}\}_{n\in\mathbb{N}})\oplus
\underline{\mbox{LIMIT}}(\{y^{(n)}\}_{n\in\mathbb{N}})\preceq
x^{(N)}+y^{(N)}+\frac{2}{10^k}.\end{align} Finally a standard
application of Lemma \ref{eqvi test improved} to (\ref{73}) and
(\ref{74}) gives the desired result, we are done.
\end{proof}

\begin{theorem}\label{cauchymultiplicative}
For any two Cauchy sequences of rational numbers
$\{x^{(n)}\}_{n\in\mathbb{N}}$, $\{y^{(n)}\}_{n\in\mathbb{N}}$,
\[\kappa\big(\{x^{(n)}\}_{n\in\mathbb{N}}\widehat{\boxtimes}\
\{y^{(n)}\}_{n\in\mathbb{N}}\big)=
\kappa\big(\{x^{(n)}\}_{n\in\mathbb{N}}\big)\widehat{\otimes}\
\kappa\big(\{y^{(n)}\}_{n\in\mathbb{N}}\big).\]
\end{theorem}

\begin{proof}

By Formula (13) in Subsection \ref{subsection Golden Flowers} and
Theorem \ref{aaaaaaaa}, to prove this theorem we may assume both
$\underline{\mbox{LIMIT}}(\{x^{(n)}\}_{n\in\mathbb{N}})\succeq0.000000\cdots$
and
$\underline{\mbox{LIMIT}}(\{y^{(n)}\}_{n\in\mathbb{N}})\succeq0.000000\cdots$,
which easily implies that $\{x^{(n)}\}_{n\in\mathbb{N}}$ and
$\{y^{(n)}\}_{n\in\mathbb{N}}$ are eventually positive, that is,
there is an starting index $M\in\mathbb{N}$ such that $\forall n\geq
M$, $x^{(n)}\succeq0.000000\cdots$, $y^{(n)}\succeq0.000000\cdots$.

Similar to the proof of the previous theorem,  for arbitrary
$k\in\mathbb{N}$ one can find a starting index $N\geq M$ such that
$\forall m,n \geq N$, $|x^{(m)}-x^{(n)}|<\frac{1}{10^k}$,
$|y^{(m)}-y^{(n)}|<\frac{1}{10^k}$, and also
\begin{align}
x^{(N)}-\frac{1}{10^k}&\preceq \underline{\mbox{LIMIT}}(\{x^{(n)}\}_{n\in\mathbb{N}})\preceq x^{(N)}+\frac{1}{10^k},\label{75}\\
y^{(N)}-\frac{1}{10^k}&\preceq \underline{\mbox{LIMIT}}(\{y^{(n)}\}_{n\in\mathbb{N}})\preceq y^{(N)}+\frac{1}{10^k},\label{76}\\
x^{(N)}\cdot y^{(N)}-\frac{2z}{10^k}-\frac{1}{100^k}&\preceq
\underline{\mbox{LIMIT}}(\{x^{(n)}\cdot
y^{(n)}\}_{n\in\mathbb{N}})\preceq x^{(N)}\cdot
y^{(N)}+\frac{2z}{10^k}+\frac{1}{100^k},\label{77}
\end{align}
where $z\in\mathbb{N}$ is any common upper bound for the sequences
$\{|x^{(n)}|\}_{n\in\mathbb{N}}$, $\{|y^{(n)}|\}_{n\in\mathbb{N}}$.

\textsc{Case 1}: Suppose $x^{(N)}-\frac{1}{10^k}\geq0$ and
$y^{(N)}-\frac{1}{10^k}\geq0$. Then applying Theorem \ref{keep order
2} to (\ref{75}) and (\ref{76}) yields
\begin{align}\label{78}s-\frac{2z}{10^k}-\frac{1}{100^k}\preceq
\underline{\mbox{LIMIT}}(\{x^{(n)}\}_{n\in\mathbb{N}})\otimes
\underline{\mbox{LIMIT}}(\{y^{(n)}\}_{n\in\mathbb{N}})\preceq
s+\frac{2z}{10^k}+\frac{1}{100^k},\end{align} where $s$ stands for
$x^{(N)}\cdot y^{(N)}$ purely for the sake of simplicity. Next a
standard application of Lemma \ref{eqvi test improved} to (\ref{77})
and (\ref{78}) gives the desired result.

\textsc{Case 2}: Suppose $x^{(N)}-\frac{1}{10^k}<0$. Now we revise
(\ref{75})$\sim$(\ref{77}) slightly to
\begin{align}
0&\preceq \underline{\mbox{LIMIT}}(\{x^{(n)}\}_{n\in\mathbb{N}})\preceq \frac{2}{10^k},\label{79}\\
0&\preceq
\underline{\mbox{LIMIT}}(\{y^{(n)}\}_{n\in\mathbb{N}})\preceq
z+\frac{1}{10^k},\label{710}\\
0&\preceq \underline{\mbox{LIMIT}}(\{x^{(n)}\cdot
y^{(n)}\}_{n\in\mathbb{N}})\preceq
\frac{3z}{10^k}+\frac{1}{100^k}.\label{711}
\end{align}
Then applying Theorem \ref{keep order 2} to (\ref{79}) and
(\ref{710}) yields
\begin{align}\label{712}0\preceq
\underline{\mbox{LIMIT}}(\{x^{(n)}\}_{n\in\mathbb{N}})\otimes
\underline{\mbox{LIMIT}}(\{y^{(n)}\}_{n\in\mathbb{N}})\preceq
\frac{2z}{10^k}+\frac{2}{100^k},\end{align} Finally a standard
application of Lemma \ref{eqvi test improved} to (\ref{711}) and
(\ref{712}) gives the desired result.

\textsc{Case 3}: Suppose $y^{(N)}-\frac{1}{10^k}<0$. The proof is
fully identical to that in Case 2.

This concludes the whole proof of the theorem.
\end{proof}

Given two Cauchy sequences of rational numbers
$\{x^{(n)}\}_{n\in\mathbb{N}}$, $\{y^{(n)}\}_{n\in\mathbb{N}}$, we
say $\{x^{(n)}\}_{n\in\mathbb{N}}$ is less than or equal to
$\{y^{(n)}\}_{n\in\mathbb{N}}$, denoted by
$\{x^{(n)}\}_{n\in\mathbb{N}}\widehat{\lesssim}\{y^{(n)}\}_{n\in\mathbb{N}}$,
if either
$\{x^{(n)}\}_{n\in\mathbb{N}}\approx\{y^{(n)}\}_{n\in\mathbb{N}}$ or
the sequence $\{y^{(n)}-x^{(n)}\}_{n\in\mathbb{N}}$ is eventually
positive, that is, there exists an $n\in\mathbb{N}$ such that
$\forall n\geq N$, $y^{(n)}-x^{(n)}>0$. It is also very easy to
verify that
$$\{x^{(n)}\}_{n\in\mathbb{N}}\ \widehat{\lesssim}\ \{y^{(n)}\}_{n\in\mathbb{N}}\Leftrightarrow
\kappa(\{x^{(n)}\}_{n\in\mathbb{N}})\ \widehat{\preceq}\
\kappa(\{y^{(n)}\}_{n\in\mathbb{N}}).$$ Thus with Theorem \ref{keep
order 33}, Theorem \ref{keep order 333333}, Theorem \ref{cauchy72},
Theorem \ref{cauchyadditive} and Theorem \ref{cauchymultiplicative},
$(\mathbb{R},\widehat{\oplus},\widehat{\otimes},\widehat{\preceq})$
is isomorphic to
$(\widehat{\mathbf{CR}},\widehat{\boxplus},\widehat{\boxtimes},\widehat{\lesssim})$.

We note a noteworthy difference. When  constructing the real numbers
via decimals as before, we introduced a total order as earlier as
possible, then based on derived operations of this order, we
introduced additive and multiplicative operations. But when people
work on the Cauchy sequence approach, total order is not so
important a concept as we regard. What they prefer most is another
fundamental concept called ``\textsf{distance}". It is amazing to
see different approaches lead to the same structure, mathematics is
so harmonious!

\section{Rationals v.s. Irrationals}

According to Lemma \ref{Interesting}, there are generally three
types of elements in $\mathcal R$, that is,
\begin{itemize}
\item Type 1: $\mathbb{Q}$
\item Type 2: ${\mathcal R}_9$
\item Type 3: The others remained
\end{itemize}
According to Lemma \ref{R9Qf}, every element of second type will be
identified with an element of first type in $\mathbb{R}$. Thus no
matter in $\mathcal R$ or in $\mathbb{R}$, we should pay attention
to the third type of elements, which has deserved not so better a
name in history.

\begin{defn}
An element of $\mathcal R$ is called  \textsf{irrational} if it does
not belong to the first two types of elements.
\end{defn}

\begin{theorem} An element
$x_0.x_1x_2x_3\cdots\in{\mathcal R}$ is irrational if and only if it
cannot end with infinite recurrence of a block of digits.
\end{theorem}

\begin{proof}
By the pigeonhole principle, it is easy to prove that any rational
number $\alpha=\frac{q}{p}$ with $p\in\mathbb{N}$ and
$q\in\mathbb{Z}$ must end with infinite recurrence of a block of
length $\leq p$. Obviously, any element of ${\mathcal R}_9$ is also
of such property. Thus to prove this theorem, it suffices to show
any element ends with infinite recurrence of a block of digits must
belong to the first two types. To this aim, without loss of
generality let us suppose
\[x=x_0.x_1\cdots x_k\underbrace{y_1\cdots y_s}\underbrace{y_1\cdots y_s}\underbrace{y_1\cdots y_s}\underbrace{y_1\cdots y_s}
\underbrace{y_1\cdots y_s}\underbrace{y_1\cdots
y_s}\cdots\in{\mathcal R}\] with $y_s<9$ (please think why we can
impose this condition), then we need only prove $x$ belongs to the
first type, that is, $x\in\mathbb{Q}$. Now we define a rational
number
\[z\triangleq [x]_k+\frac{\displaystyle 1}{\displaystyle 10^k}\cdot\frac{\displaystyle \sum_{j=1}^{s}y_j\cdot10^{s-j}}{\displaystyle 10^s-1}
=\frac{\displaystyle(\sum_{i=0}^{k}x_i\cdot10^{k-i})\cdot(10^s-1)+\sum_{j=1}^{s}y_j\cdot10^{s-j}}{\displaystyle
10^k\cdot(10^s-1)}.\] Then in the following we devise two sequences
of rational numbers with left hand $\{\mu_n\}_{n\in\mathbb{N}}$, and
right hand $\{\omega_n\}_{n\in\mathbb{N}}$ satisfying the desired
inequalities\footnote{Let us see from a special example to know why
these inequalities are true, from which one can easily prove the
general cases. For example, suppose $x=0.232323\cdots$ and
$z=\frac{23}{99}$.  Take the third of this sequence of inequalities
for example, multiplying each parts by$\frac{99}{23}$ gives
$$0.999999\leq1\leq0.999999+0.000001\cdot\frac{99}{23},$$ which is
obviously true.}
\begin{align*}
x_0.x_1\cdots x_k\underbrace{y_1\cdots y_s}&\leq z\leq
x_0.x_1\cdots x_k\overbrace{y_1\cdots y_{s-1}(y_s+1)}\\
x_0.x_1\cdots x_k\underbrace{y_1\cdots y_s}\underbrace{y_1\cdots
y_s}&\leq z\leq x_0.x_1\cdots x_k\underbrace{y_1\cdots y_s}\overbrace{y_1\cdots
y_{s-1}(y_s+1)}\\
 x_0.x_1\cdots x_k\underbrace{y_1\cdots
y_s}\underbrace{y_1\cdots y_s}\underbrace{y_1\cdots y_s}&\leq z\leq
x_0.x_1\cdots x_k\underbrace{y_1\cdots y_s}\underbrace{y_1\cdots
y_s}\overbrace{y_1\cdots y_{s-1}(y_s+1)}\\
\cdots &\leq z \leq \cdots.
\end{align*}
 Note
$\{\mu_n\}_{n\in\mathbb{N}}$ is monotonically increasing with limit
$x$, while $\{\omega_n\}_{n\in\mathbb{N}}$ is monotonically
decreasing with  limit also $x$, thus by Lemma \ref{lemma order} we
must have $x=z$. This proves $x$ is a rational number, we are done.
\end{proof}

\section{Algebraic abstraction: Complete ordered field}

After detailed studies of three explicit approaches to the real
number system, let us explain how can we come to its algebraic
characterization by the so-called concept of \textsf{complete
ordered field}.

 If there exists only one
\textsf{real number system} in the mathematical world, then it
should be a system $(\mathbf{R},+,\times,\leq)$ with  three binary
operations $+,\times,\leq$ endowed on the same ambient set
$\mathbf{R}$ satisfying at least a few conditions:

\begin{itemize}
\item \textsf{$\mathbf{R}$ is  infinite (reason: should contain $\mathbb{Q}$ as a
subset)}
\item \textsf{$(\mathbf{R},+,\times)$ is a field (reason: $(\mathbb{R},\widehat{\oplus},\widehat{\otimes})$
is a field)}
\item \textsf{$\leq$ is a total order on $\mathbf{R}$
 and $(\mathbf{R},\leq)$ has both the least upper bound property and the
greatest lower bound property (reason:
$(\mathbb{R},\widehat{\preceq})$ is of such properties)}
\end{itemize}

An algebraic system satisfies the above conditions  can not
characterize the real number system yet. The main reason is we can
not sense any difference between the \textsf{positive} and the
\textsf{negative} parts of our desired real number system, so we
need more observation. To mention a few: addition of elements of
$\widehat{\mathbf{CR}}$ preserve order $\widehat{\lesssim}$, so are
multiplication by positive elements of $\widehat{\mathbf{CR}}$.

Is an algebraic system $(\mathbf{R},+,\times,\leq)$ satisfies the
following conditions isomorphic to those systems of Dedekind and
M\'{e}ray-Cantor-Heine?
\begin{itemize}
\item \textsf{$(\mathbf{R},+,\times)$ is an infinite field}
\item \textsf{$\leq$ is a total order on $\mathbf{R}$
 and $(\mathbf{R},\leq)$ has both the least upper bound property and the
greatest lower bound property }
\item \textsf{$x<y,z\in\mathbf{R}\Rightarrow x+z<y+z$}
\item \textsf{$x<y, z>0 \Rightarrow xz<yz$}
\end{itemize}
Yes, it is! Such a system is called a \textsf{complete ordered
field}. It does not so matter whether some of the conditions are
superfluous  or not, at least we derived an algebraic
characterization of the well-known real number system. How to prove
it? Decimal representations! (Imagine it first, prove it second!)

(Remark: The infinite cardinality assumption follows from the third
assumption and one item of the field axioms, that is, the additive
unit is not equivalent to the multiplicative unit, while the
greatest lower bound property follows from the least lower bound
property and the other three assumptions, so they are unnecessary.)

\section{Epilog: Analysis v.s. Algebra or Analysis plus Algebra ?}

When to teach mathematical majors Mathematical Analysis, instructors
need to make a serious choice.

\underline{\textsc{Dedekind cuts}}

``\textsf{Although the definition is geometrically motivated, and is
based on a `natural idea', the definition is actually far too
abstract to be efficient at the beginning of a calculus course.}"
(\cite{Leviatan}) We also don't prefer this approach. The main
reason is, once the real number system is established through this
approach, we can find almost no use of the language of Dedekind cuts
in the subsequent study of the Mathematical Analysis course.

\underline{\textsc{Cauchy sequences}}

We believe this approach is  too abstract to be acceptable for
beginners, who have not received adequate mathematical training yet.
We also  have a sense that when to construct the real number system
explicitly, people would prefer more the Dedekind cut approach than
the abstract Cauchy sequence one. Introducing it in the Functional
Analysis course would be a better choice.

\underline{\textsc{Axiomatic definition}}

As said in the Introduction, we found  many authors liked to choose,
say for example the algebraic-axiomatic definition,  then proved
that every real number has a suitable decimal representation. But
how can we persuade students believe that a real number is an
element of an algebraic system satisfying a dozen or so properties
at the beginning of a Mathematical Analysis course?

\underline{\textsc{Suggestion: Decimal (First) plus Axiomatic
(Second) approaches}}

\textsf{This paper  provides a complete approach to the real number
system via rather old decimal representations, which perfectly
matches what we have learnt a number is in high school. Once the
real number system is established through this approach, the least
upper bound property, the greatest lower bound property, the
supremum, infimum, upper limit and lower limit operations are
immediately, also naturally obtained. This shows huge advantages
over the Dedekind cut and the Cauchy sequence approaches. It would
be great if after the decimal approach, one can give an algebraic
characterization of the real number system, which help students
grasp the fundamental properties of the real numbers instead of
heavily relying on the decimal representations, also help them
understand there is only one real number system in the mathematical
world when they consult other Mathematical Analysis books.}

Although we have not reviewed many other construction of the real
number system, we show full respect to any such kind of work. We
conclude this paper with a saying by Jonathan Bennett: ``\textsf{You
can have lots of good ideas while walking in a straight line.}"

\end{document}